\documentclass{article}

\usepackage{latexsym,amssymb}
\title{\bf  Elements of the Green Ring related to Alperin's Weight Conjecture}
\author{Geoffrey R. Robinson\\
Institute of Mathematics\\
University of Aberdeen\\
Fraser Noble Building\\
Aberdeen\\
AB24 3UE\\
United Kingdom\\}

\begin{document}

\maketitle

\medskip
\begin{center}
\emph{Dedicated to the memory of Jonathan L. Alperin}
\end{center}

\begin{abstract}
In the article [1], J. L. Alperin introduced his celebrated weight conjecture, a major question in the modular representation theory of finite groups, which has now been studied for several decades. Alperin viewed this conjecture as the beginning of a ``Lie theory for the general finite group", as the title of the article indicates. This paper is intended as a tribute to Alperin's insight. In this note, we place an emphasis on projective modules
and putting the numerical questions arising from the weight conjecture into a ``module form" within the appropriate Green ring of $RG$, when $G$ is a finite group, and $R$ is a suitable complete discrete valuation ring with residue field of prime characteristic $p$. We believe that this work moves in the direction of a structural approach to Alperin's weight conjecture. 
\end{abstract}

\section{Introduction}

\medskip
The paper [8] of R. Kn\"orr and the present author is our starting point, and we develop further some of the theory begun in our earlier papers [15] and [16], using virtual projective modules in the Green ring to study Alperin's Weight Conjecture. 

\medskip
Our earlier work [16] introduced a virtual projective module $P_{1,p,G}$, defined in purely group-theoretic terms, and placed emphasis on the`` truncated conjugation module", a ``projective version" of the group algebra under the conjugation action. By contrast, in the much earlier paper [15], we had decomposed the truncated conjugation module (in fact for a block) in representation-theoretic terms via the Lefschetz conjugation module (which had its genesis in the paper [8]) by a (non-conjectural) module-theoretic formula in the appropriate Green ring which also gave a (virtual) module-theoretic encoding of Alperin's weight conjecture. Understanding the relationship  between these two decompositions is a key theme of this paper. We attempt to illustrate here this relationship, and the interaction of these decompositions with the standard theory, in some familiar examples.

\medskip
We fix a finite group $G$ and a prime $p$,  and we choose a $p$-modular system $(R,\mathbb{K},\mathbb{F})$ for $G$ which is splitting for $G$ and its subgroups. We denote the unique maximal ideal of $R$ by $J(R).$ A block $B$ of the group ring $RG$ is a minimal ($R$-free) two-sided ideal of $RG$. Blocks of $B$ are in bijection with the blocks of the group algebra $\mathbb{F}G$, realised in one direction by reduction (${\rm mod}$  $J(R)$). If we need to distinguish between the block $B$ of $RG$ and its image block of $\mathbb{F}G$, we will denote the latter by $\tilde{B}.$

\medskip
We let $k(G)$ denote the number of conjugacy classes of $G$, we let $\ell(G)$ denote the number of conjugacy classes of $p$-regular elements of $G$, and we let $\ell_{0}(G)$ denote the number of conjugacy classes of elements $y \in G$ such that $|C_{G}(y)|$ is not divisible by $p$. The number $\ell_{0}(G)$ is sometimes referred to as the number of ($p$)-defect zero classes of $G$. It is also the dimension of ${\rm Tr}_{1}^{G}(\mathbb{F}G),$ as well as being the $p$-rank (or number of determinantal divisors $1$) of the Cartan matrix for $\mathbb{F}G$.The number of isomorphism types of simple $\mathbb{F}G$-modules is also equal to $\ell(G),$ and when $B$ is a sum of blocks of $RG$, we let $\ell(B)$ denote the number of isomorphism types of simple $B$-modules.

\medskip
For $B$ a block of $RG$, we let $\ell_{0}(B)$ denote the $\mathbb{F}$-dimension of ${\rm Tr}_{1}^{G}(\tilde{B}),$ which is also the $p$-rank of the Cartan matrix of $B$.

\medskip
We let $f_{0}(G)$ denote the number of blocks of defect zero of $RG$, which is also the number of isomorphism types of projective simple $\mathbb{F}G$-modules. We let $\mathcal{Q}$ denote a set of representatives for the conjugacy classes of $p$-subgroups of $G$ (including the trivial subgroup). If $B$ is a sum of blocks of $RG$, and $Q \in \mathcal{Q},$ we let $B(Q)$ denote the sum of those blocks of $RN_{G}(Q)/Q$ whose modules lie in Brauer correspondents of $B$ for $N_{G}(Q)$ when inflated to $RN_{G}(Q)$-modules. We let $f_{0}(B(Q))$ denote the number of block summands of $B(Q)$ which have defect zero.

\medskip 
We will deal mostly with (virtual) projective $RG$-modules.  One reason for this choice is that the ideal of the character ring of $G$ consisting of virtual characters which vanish on all $p$-singular elements of $G$ faithfully represents the ideal of the Green ring of $RG$ consisting of $\mathbb{Z}$-linear combinations of projective $RG$-modules. This latter ideal is in turn isomorphic (as Abelian group) to the ideal of the Green ring of $\mathbb{F}G$ consisting of $\mathbb{Z}$-linear combinations of projective $\mathbb{F}G$-modules. 

\medskip
The fact that the ideal of the Green ring of $\mathbb{F}G$-modules  consisting of virtual projective modules is invariant under multiplication by simple $\mathbb{F}G$-modules in a ``projective module preserving fashion" lifts to a useful property of the ideal of the Green ring of $RG$ consisting of projective $RG$-modules: characters of projective $RG$-modules multiplied by Brauer characters of simple $\mathbb{F}G$-module are characters of projective $RG$-modules (note that such products make sense since characters of projective $RG$-modules vanish on all $p$-singular elements).

\medskip
Another advantage of this viewpoint is that when $G$ has a non-trivial normal $p$-subgroup $Q$, there is an isomorphism of Abelian groups between the $\mathbb{Z}$ -span (in the Green ring of $RG$) of the projective indecomposable $RG$-modules, and the $\mathbb{Z}$-span (in the Green ring of $RG/Q$) of the projective indecomposable $RG/Q$-modules, which (crucially for us), is module preserving in both directions: if $M$ is a projective $RG$-module, then $M^{Q}$ is the corresponding projective $RG/Q$-module under this isomorphism, and if $N$ is a projective $RG/Q$-module affording character $\theta$, then the corresponding $RG$-module $\tilde{N}$ under this isomorphism is the projective cover of $N$ (when $N$ is inflated to an $RG$-module). 

\medskip
In general, if a virtual projective $RG/Q$-module affords virtual character $\theta,$ then the corresponding virtual projective $RG$-module is the unique virtual projective $RG$-module which affords virtual character agreeing with $\theta\gamma_{Q}$ on $p$-regular elements (and necessarily vanishing on all $p$-singular elements), where $\gamma_{Q}$ is the permutation character of $G$ afforded by its conjugation action on $Q$. This bijection between virtual projectives of $RG$ and virtual projectives of $RG/Q$ is due to W.F. Reynolds [11].

\medskip
This isomorphism allows us to speak in a well-defined way of the (virtual) projective cover of a virtual projective $RG/Q$-module $X$. If $X$ is a virtual projective $RG/Q$-module affording virtual character $\chi$, then we will refer to the unique virtual projective $RG$-module which affords virtual character vanishing on all $p$-singular elements and agreeing with $\chi \gamma_{Q}$ on all $p$-regular elements as the (virtual) projective cover of $X$ as (virtual)  projective $RG$-module, and denote it by $P[X].$

\medskip
We recall that $\ell(G)$ is also the rank of the $\mathbb{Z}$-submodule (and ideal) of the character ring of $G$ consisting of virtual characters vanishing on all $p$-singular elements.

\medskip
\noindent  {\bf Remark 1.1:} We state Alperin's weight conjecture and its blockwise refinement.

\medskip
\noindent {\bf Alperin's Weight Conjecture:} \emph{We have $\ell(G) = \sum_{Q \in \mathcal{Q}} f_{0}(N_{G}(Q)/Q).$}

\medskip
\noindent {\bf Blockwise Refinement:}  \emph{For each block $B$ of $RG$,  we have $\ell(B) = \sum_{Q \in \mathcal{Q} } f_{0}(B(Q)).$}

\medskip
\noindent {\bf Definition 1.2:} The \emph{Truncated Conjugation Module for G} is the $RG$-module, denoted by $T_{c,p,G}$, which is the lift to an $RG$-module of the projective $\mathbb{F}G$-module $\oplus_{S} {\rm Hom}_{\mathbb{F}}(S,P(S))$, where $S$ runs over the isomorphism types of simple $\mathbb{F}G$-modules, and $P(S)$ denotes the projective cover of $S$.

\medskip
The projective $RG$-module $T_{c,p,G}$ affords a character $\Lambda_{c,p,G},$  which vanishes on all $p$-singular elements of $G$, and takes value 
$|C_{G}(x)|$ at each $p$-regular $x \in G.$ This can be seen by block orthogonality relations due to Brauer, or by considering the diagonal matrix $\overline{\Phi}^{t}C\Phi,$ where $C$ is the Cartan matrix and $\Phi$ is the Brauer character table of $G$ (in matrix form). 

\medskip
We note that $$T_{c,p,G} = P[T_{c,p,G/Q}]$$ whenever $Q$ is a normal $p$-subgroup of $G$, as may be seen by  consideration of the characters they afford.

\medskip
Note that the multiplicity of the projective cover of the trivial $RG$-module as a summand of ${\rm T}_{c,p,G}$ is the number of isomorphism types of simple $\mathbb{F}G$-modules, which is $\ell(G).$ This can also be seen from the character $\Lambda_{c,p,G}.$

\medskip
We consider the truncated conjugation module as a ``projectivisation" of the $RG$-module afforded by the conjugation action of $G$ on $RG$.  

\medskip
\noindent {\bf Remark 1.3:} It may be instructive to note that  $T_{c,p,G}$ is the regular module of $RG$ if and only if $G$ is a $p$-group, and that $T_{c,p,G}$ is isomorphic to $RG$ viewed as a module under the conjugation action of $G$ if and only if $G$ is a $p^{\prime}$-group.

\medskip
In fact, in [15], we had already considered the truncated conjugation module of a block $B$ of $\mathbb{F}G$, which is the projective $\mathbb{F}G$-module $$\oplus_{S} {\rm Hom}_{\mathbb{F}}(S,P(S)),$$ where $S$ runs over the isomorphism types of simple (over $\mathbb{F}$) $B$-modules, though we did not use this terminology there. This projective $\mathbb{F}G$-module lifts uniquely to a projective ($R$-free) $RG$-module, which we will here denote by $T_{c,p,B}$ (and which we will still refer to as the truncated conjugation module of $B$).

\medskip
\noindent {\bf Remark 1.4:} $T_{c,p,B}$ affords character $\Lambda_{c,p,B}$ which vanishes on all $p$-singular elements of $G$, and takes value
 $\sum_{ \chi \in {\rm Irr}(B)}|\chi(x)|^{2}$ for each $p$-regular $x \in G$ (this also follows by block orthogonality relations and by properties of the Cartan matrix).  The multiplicity of the projective cover of the trivial $RG$-module as a summand of ${\rm T}_{c,p,B}$ is clearly $\ell(B).$

\medskip
Note that the indecomposable summands of ${\rm T}_{c,p,B}$ do not in general all lie in the block $B$ (the projective cover of the trivial module will always occur as a summand of $T_{c,p,B}$, for example, while $B$ need not be the principal block).

\medskip
\noindent {\bf Remark 1.5:} In [8], it was proved that  (in the Green ring for $\mathbb{F}G$), a certain alternating sum of induced modules is a virtual projective module, and that Alperin's weight conjecture is equivalent to proving that the number of ``fixed points" (strictly, the number of occurrences of the projective cover of the trivial $RG$-module) on this virtual projective module is $f_{0}(G)$. The blockwise version says that the blockwise version of Alperin's weight conjecture is equivalent to proving that the number of ``fixed-points" on the block version of this virtual module is zero for blocks of positive defect, and is one  for blocks of $\mathbb{F}G$ of defect zero (the second of these conjectured equalities is almost immediate, it is the first that requires proof).

\medskip
We leave this virtual module undefined for the moment, until we have discussed the $p$-subgroup complex of $G$. In [15], we referred to this virtual projective module as the \emph{Lefschetz Conjugation Module of $G$}, and also defined the Lefschetz conjugation module of a block $B$ of $\mathbb{F}G$. In the present work, we prefer to deal with the lifts of these virtual projective modules to virtual projective $RG$-modules. The notation we will adopt here is that the Lefschetz conjugation module of a direct sum $B$ of blocks  of $RG$ is denoted by $\mathcal{L}_{c}(B).$ We had previously (in [8]) defined a block $B$ of $RG$ of positive defect to be of $\emph{Lefschetz type}$ if the coefficient of the projective cover of the trivial module in $\mathcal{L}_{c}(B)$(when expressed as a $\mathbb{Z}$-combination  of projective $RG$-modules) is zero. This allows us to state succinctly that the blockwise version of  Alperin's weight conjecture is true for the prime $p$ if and only if, whenever $G$ is a finite group, every block of positive defect of $RG$ is of Lefschetz type.

\medskip
It should be noted, however, that a priori, a particular block of Lefschetz type need not satisfy the Alperin weight conjecture, but if it does not, then there is a a block of some proper section of $G$ which is not of Lefschetz type.

\medskip
Motivated by the results of [8], in [15] we used a M\"obius inversion-type argument to express the projective $\mathbb{F}G$-module
$$\oplus_{S}{\rm Hom}_{\mathbb{F}}(S,P(S))$$ as a sum of virtual projective $\mathbb{F}G$-modules, each induced from virtual projective modules 
of $p$-local subgroups of $G$, in a way which reflects in the Green ring the numerical equality predicted by Alperin's Weight Conjecture. However, rather than being conjectural, this is a genuine equality in the Green ring of $\mathbb{F}G$.

\medskip
With the notation and terminology we have adopted here, this formula may be written (when lifted to $RG$-modules) as follows:

\medskip
\noindent {\bf Remark 1.6:} For each block $B$ of $RG$, we have $$T_{c,p,B} = \sum_{Q \in \mathcal{Q}} {\rm Ind}_{N_{G}(Q)}^{G}(P[\mathcal{L}_{c}(B(Q))]).$$ 

\medskip
This may be compared with the form of the blockwise version of Alperin's weight conjecture given in 1.1, bearing in mind that the multiplicity of the projective cover of the trivial $RG$-module on the left-hand side is $\ell(B).$

\medskip
Once this equality in the Green ring is established, it is easy to read off the fact that the blockwise form of Alperin's Weight Conjecture is equivalent to proving that the number of occurrences of the projective cover of the trivial module in the virtual projective module $\mathcal{L}_{c}(B(Q))$ is (always) $f_{0}(B(Q))$. Note that if $B$ itself is not of defect zero, then this prediction should imply that all occurrences of the projective cover of the trivial $RG$-module on the right side come from non-trivial $p$-subgroups $Q$, which means that $\ell(B)$ is (strictly) $p$-locally determined. 

\medskip
We regard the expression in 1.6 as a \emph{representation-theoretic decomposition of $T_{c,p,B}$}, which may be summed over all blocks of $RG$ to obtain a representation theoretic decomposition of $T_{c,p,G}.$

\medskip
By contrast, in [16], we decomposed the virtual projective module $T_{c,p,G}$ in terms of purely group-theoretic information (it is not clear to us at present whether there might be a block-theoretic refinement of such a decomposition).

 \medskip
\noindent {\bf Remark 1.7: } In [16], we expressed the character $\Lambda_{c,p,G}$ in the form $$\Lambda_{c,p,G} = \sum_{y \in G_{p^{\prime}}/G} {\rm Ind}_{C_{G}(y)}^{G}( \Psi_{1,p,C_{G}(y)}),$$ where for any finite group $H$,  the generalized character $\Psi_{1,p,H}$ vanishes on all $p$-singular elements of $H$, and takes value $|H_{p} \cap C_{H}(h)|$ for each $p$-regular $h \in H$ (here, $H_{p}$ denotes the set of $p$-elements of $H$ and $G_{p^{\prime}}$ denotes the set of $p$-regular elements of $G$).  

\medskip
We also conjectured in [16] that the unique virtual projective $RH$-module $P_{1,p,H}$ affording virtual character $\Psi_{1,p,H}$  is in fact a projective $RH$-module, and proved the conjecture  in some cases. We consider the above expression for the character $\Lambda_{c,p,G}$  as a \emph{group-theoretic} description (or, equivalently, as a group-theoretic description of the projective module ${\rm T}_{c,p,G}$), since the generalized characters $\Psi_{1,p,H}$ (and hence the virtually projective $RH$-modules $P_{1,p,H}$) are described in purely group-theoretic terms. Notice however that the original definition of the module $T_{c,p,G}$ is intrinsically representation-theoretic.

\medskip
In [16], we also gave an explicit expression for $P_{1,p,G}$ as  a sum of virtual projective modules in the Green ring of $RG$.  Before we can describe this, we need to discuss $p$-subgroup complexes, and the associated (virtual) Steinberg modules which arise. 

\medskip
\noindent {\bf Remark 1.8: } For the sake of economy of notation, when $y$ is a $p$-regular element of $G$, we denote the induced generalized character 
$${\rm Ind}_{C_{G}(y)}^{G}( \Psi_{1,p,C_{G}(y)})$$ by $\Psi_{y,p,G},$ and the associated virtual projective $RG$-module by $P_{y,p,G}.$ It is easy to check that $P_{y,p,G}$ is the inflation to an $RG$-module of $P_{1,p,G/\langle y \rangle}$ when $y$ is a central $p$-regular element of $G$.

\medskip
\section{The relationship to the $p$-subgroup complex and the generalized Steinberg (virtual) module of P.J. Webb}

\medskip
Let $\mathcal{S}_{p}(G)$ denote the $p$-subgroup complex of the finite group $G$ (we include the empty simplex, though it is sometimes convenient to identify it with the singleton $1_{G}$, with chain stabilizer $G$). A non-empty simplex in $\mathcal{S}_{p}(G)$ is a strictly increasing chain $$\sigma = Q_{1} < Q_{2} < \ldots < Q_{n}$$ of non-trivial $p$-subgroups $Q_{i}$ of $G$, and we use $|\sigma|$ to denote the number of non-trivial $p$-subgroups in $\sigma$. Clearly, $G$  acts by conjugation on $\mathcal{S}_{p}(G)$, and we let $G_{\sigma}$ denote the stabilizer of $\sigma$ under this action.

\medskip
The \emph{Steinberg virtual module}, denoted ${\rm St}_{p}(G)$, and defined by P.J. Webb,[17], is (apart from a different sign convention) the element  
$$\sum_{ \sigma \in \mathcal{S}_{p}(G)/G} (-1)^{|\sigma|} {\rm Ind}_{G_{\sigma}}^{G}(R)$$ of the Green ring of $RG$, (where $R$ denotes the trivial $RG$-module). Webb proved that ${\rm St}_{p}(G)$ is always a virtual projective module in the Green ring of $RG$. If $O_{p}(G) \neq 1,$ then (as noted by Webb), the $p$-subgroup complex is $G$-contractible, and ${\rm St}_{p}(G) = 0.$

\medskip
\noindent {\bf Remark 2.1: } In [16] it was proved that the generalized character $\Psi_{1,p,G}$ is afforded by the virtual module 
$$(\ast) \: P_{1,p,G} = \sum_{Q \in \mathcal{Q}} {\rm Ind}_{N_{G}(Q)}^{G}( P[{\rm St}_{p}(N_{G}(Q)/Q)]).$$ 

\medskip
We note that, since $${\rm St}_{p}(N_{G}(Q)/Q)=0$$ if $O_{p}(N_{G}(Q)/Q) \neq 1,$ we may restrict the subgroups $Q$ considered to those which satisfy $Q = O_{p}(N_{G}(Q)),$ the so-called $p$-\emph{radical} subgroups of $G$. We let $\mathcal{R}$ denote a set of representatives for the conjugacy classes of $p$-subgroups of $G$ such that $Q = O_{p}(N_{G}(Q)).$

\medskip
\noindent {\bf Remark 2.2:} We note that ${\rm St}_{p}(G)$ is the trivial module $R$ if $G$ is a $p^{\prime}$-group, and that the multiplicity of the projective cover of the trivial module in ${\rm St}_{p}(G)$ is zero whenever $p$ divides $|G|$. Hence it follows that the multiplicity of the projective cover of the trivial module in $P_{1,p,G}$ is precisely one, coming from the case $Q \in {\rm Syl}_{p}(G)$ in the above expression ($\ast$) for $P_{1,p,G}$.

\medskip
It was also proved in [16] that we have 
$$T_{c,p,G} = \sum_{y \in G_{p^{\prime}}/G} P_{y,p,G}$$ in the Green ring of $RG$. It follows that each $p$-regular conjugacy class of $G$ contributes  one copy of the projective cover of the trivial module for $RG$ to this expression for $T_{c,p,G}$, which comes (for a class with representative $y$) from 
$${\rm Ind}_{C_{G}(y)}^{G}(P[{\rm St}_{p} ((C_{G}(y) \cap N_{G}(Q))/Q)]),$$ where $Q \in {\rm Syl}_{p}(C_{G}(y)).$

\medskip
\noindent {\bf Remark 2.3: } Using the $p$-subgroup complex, we can now state that the Lefschetz conjugation module $\mathcal{L}_{c}(RG)$ is the element $$\sum_{ \sigma \in \mathcal{S}_{p}(G)/G} (-1)^{|\sigma|} {\rm Ind}_{G_{\sigma}}^{G}(RG_{\sigma})$$ of the Green ring of $RG$,  where $RG_{\sigma}$ is viewed as an $RG_{\sigma}$ module under conjugation action for each $\sigma.$ By the results of [8], (lifted to $RG$-modules), $\mathcal{L}_{c}(RG)$ is a virtual projective $RG$-module, and it is proved in [8] (in fact block by block) that the multiplicity of the projective cover of the trivial module in $\mathcal{L}_{c}(RG)$ is equal to  $$\sum_{\sigma \in \mathcal{S}_{p}(G)/G} (-1)^{|\sigma|} k(G_{\sigma}),$$ where $k(G_{\sigma})$ is the number of conjugacy classes of $G_{\sigma}$, and that $k(G_{\sigma})$ can be replaced by $\ell(G_{\sigma})$ throughout (or by $\ell_{0}(G_{\sigma})$ throughout). Furthermore, it is proved there that Alperin's weight conjecture (in non-blockwise form), is equivalent to the statement that this number is equal (for every finite group $G$) to $f_{0}(G)$, and that a similar statement holds block by block.

\medskip
Since virtual projective modules of the Green ring of $RG$ are uniquely specified by their virtual characters, and since, for a block $B$ of $RG$, the genuine  $RG$-modules $B$ (with conjugation action) and the projective module $T_{c,p,B}$ afford characters which agree on $p$-regular elements, we have:

\medskip
\noindent {\bf Lemma 2.4:} \emph{ For each block $B$ of $RG$, we have
 $$\mathcal{L}_{c}(B) = \sum_{\sigma \in \mathcal{S}_{p}(G)/G} (-1)^{|\sigma|}{\rm Ind}_{G_{\sigma}}^{G}(T_{c,p,B_{\sigma}}),$$ where $B_{\sigma}$ denotes the sum of the Brauer correspondent blocks for $B$ for the chain stabilizer $G_{\sigma}$.}

\medskip
In this statement, we have now expressed $\mathcal{L}_{c}(B)$ as an alternating sum of (individually) projective modules. Notice also that, counting occurrences of the projective cover of the trivial module on each side shows that the number of occurrences of the projective cover of the trivial module in $\mathcal{L}_{c}(B)$ is equal to $$\sum_{\sigma \in \mathcal{S}_{p}(G)/G} (-1)^{|\sigma|} \ell(B_{\sigma}).$$ It was also proved in [8] that the the number of occurrences of the projective cover of the trivial module in $\mathcal{L}_{c}(B)$ is equal to $$\sum_{\sigma \in \mathcal{S}_{p}(G)/G} (-1)^{|\sigma|} \ell_{0}(B_{\sigma}).$$

\medskip
\noindent {\bf Remark 2.5:} While we have conjectured in [16] that the virtual projective module $P_{1,p,G}$ is always  a module, it is definitely not  the case that the virtual projective module  $\mathcal{L}_{c}(B)$ is always a module. This is evident because when $p$ divides $|G|$, the virtual module $\mathcal{L}_{c}(RG)$ has virtual $R$-rank $$|G|\left( \sum_{\sigma \in \mathcal{S}_{p}(G)/G}(-1)^{|\sigma|}\right) = 0.$$

\medskip
\noindent {\bf Remark 2.6:} Let us summarize the interdependency between the various elements of the Green ring we have defined so far:

\medskip
\noindent i) From 2.4, we may conclude that for each block (or sum of blocks) $B$ of $RG$, the virtual projective module $\mathcal{L}_{c}(B)$ is determined by $T_{c,p,B}$ + (strictly $p$-local information built from $p$-local versions of $T_{c,p,\ast}$).

\medskip
\noindent ii) On the other hand, from 1.6, we may conclude that the projective module $T_{c,p,B}$ is determined by $\mathcal{L}_{c}(B)$ + (strictly $p$-local information built from (sectional) $p$-local versions of $\mathcal{L}_{c}(\ast)$).

\medskip
\noindent iii) From 1.7, we may conclude that the projective $RG$-module $T_{c,p,RG}$ is built from modules $P_{y,p,G}$ where $y$ runs over $p$-regular conjugacy class representatives for $G$, and each $P_{y,p,G}$ has the form ${\rm Ind}_{C_{G}(y)}^{G}(P_{1,p,C_{G}(y)}).$

\medskip
\noindent iv) From 2.1, we may conclude that for any finite group $H$, the virtual projective $RH$-module $P_{1,p,H}$ is built from the (virtual projective) Steinberg module ${\rm St}_{p}(H)$+(strictly $p$-local information built from (sectional) $p$-local versions of ${\rm St}_{p}(\ast)$).

\medskip
\noindent v) In [16], it was also proved that for any finite group $H$, the virtual projective module ${\rm St}_{p}(H)$ is built from $P_{1,p,H}$ + (strictly $p$-local information built from $p$-local versions of $P_{1,p,\ast}).$

\medskip
To avoid possible ambiguity, we emphasize that when $G$ is a $p^{\prime}$-group, the virtual modules $P_{1,p,G}$ and ${\rm St}_{p}(G)$ are both equal to the trivial module $R$,  and the virtual modules $\mathcal{L}_{c}(B)$ and $T_{c,p,B}$ are both equal to $B$ as $RG$-module under conjugation action when $B$ is a sum of blocks of $RG$. 

\medskip
When $O_{p}(G) \neq 1,$ it is the case that ${\rm St}_{p}(G)$ and $\mathcal{L}_{c}(B)$ are the zero module, while $P_{1,p,G}$ is the projective cover as virtual $RG$-module of the virtual projective $RG/O_{p}(G)$-module $P_{1,p,G/O_{p}(G)}$ and also (whenever $B$ is a sum of blocks of $RG$), the projective $RG$-module $T_{c,p,B}$ is the projective cover of $T_{c,p,B(O_{p}(G))},$ where $B(O_{p}(G))$ is the sum of blocks of $G/O_{p}(G))$ whose projective indecomposable modules inflate to modules in $B$.

\medskip
Our next result addresses the question whether the virtual module $P_{1,p,G}$ may be factored in a specific manner in the Green ring.

\medskip 
\noindent {\bf Theorem 2.7 :} \emph{Suppose that we have $$P_{1,p,G} = U^{\ast} \otimes U$$ in the Green ring of $RG$,  where $U$ is a virtual projective $RG$-module. Then there is a sign $\epsilon$ such that $\epsilon U$ is the unique projective indecomposable module $V$  in a  block $B$ of  defect zero of $RG$.}

\medskip
\noindent {\bf Proof:} Since virtual projective $RG$-modules are determined up to isomorphism by their characters, $U^{\ast} \otimes U$ is a well defined virtual projective module in the Green ring of $RG$ whenever $U$ is a virtual projective module in the Green ring of $RG$. 

\medskip
Suppose that we have $P_{1,p,G} = U^{\ast} \otimes U$ for $U$ as in the statement of Theorem 2.7. Then, as $P_{1,p,G}$ contains the projective cover of the trivial module with multipliicity one when expressed as a $\mathbb{Z}$-linear combination of projective indecomposables in the Green ring of $RG$, we conclude that the virtual character afforded by the virtual module $U^{\ast} \otimes U$ contains the trivial character with multiplicity exactly one. Hence there is a sign $\epsilon$ such that $V = \epsilon U$  is a projective $RG$-module affording an irreducible character. Then $V$ is the unique projective indecomposable $RG$-module in some block $B$ of defect zero of $RG$, as claimed.

\medskip
\noindent {\bf Remark 2.8:} One motivation for stating Theorem 2.7 is that when $G$ is a finite simple group with a $(B,N)$-pair of characteristic $p$, the generalized character $\Psi_{1,p,G}$ is the square of the Steinberg character of $G$, and in the Green ring of $RG$, we have $$P_{1,p,G} = St_{p}^{2}.$$ For the general finite group, this (virtual) module factorization appears to be the exception, rather than the rule, even in cases where $\Psi_{1,p,G}$ is known to be a character afforded by a projective $RG$-module. The theorem certainly implies that there can be no such factorization if $RG$ has no block of defect zero.

\medskip
For a specific example, when $G = A_{5}$ and $p =3,$ we see that $\Psi_{1,3,G}$ takes value $1$ on all non-identity $\{2,5\}$-elements, value $0$ on all non-identity $3$-elements, and value $21$ on the identity. It is clear then that there is no virtual character $\theta$ of $G$ such that $\theta \overline{\theta} = \Psi_{1,3,G},$ as that would require $\theta(1)^{2} = 21.$ It is though the case that $\Psi_{1,3,G} = 1 + \alpha + \beta+ 2\gamma+ \chi,$ where
$\alpha$ and $\beta$ are irreducible of degree $3$,   while $\gamma$ is irreducible of degree $5$ and $\chi$ is irreducible of degree $4$, so that $\Psi_{1,3,G}$ is  a character. In fact $\Psi_{1,3,G}$ is the character afforded by $$P_{1} \oplus P_{2} \oplus P_{3} \oplus P_{4},$$ where $P_{1}$ is the projective cover of the trivial module, $P_{2}$ is the rank $9$ projective indecomposable summand of ${\rm Ind}_{N_{G}(S)}^{G}(R)$, and $P_{3}$ and $P_{4}$ are the two non-isomorphic projective indecomposable $RG$-modules of rank $3$.

\medskip
We will next explore how the (virtual) module-theoretic descriptions, and the numerical equalities they (should) encode, behave in some familiar groups.

\medskip
\section{The case of $G \cong {\rm SL}(2,2^{n}), n \geq 2$ and other examples}

\medskip
In this section, we first analyse the case when $G \cong {\rm SL}(2,2^{n})$ in the case $p =2$, and determine the virtual modules $P_{1,2,G}$ (which are in fact modules in this case. That they are modules for this $G$ follows from the results of [16], but we will see in this case that it also follows from the general theory). We work out the details over the field $\mathbb{F}$. Since we are dealing exclusively with virtual projective modules, these results lift easily to the required statements about $RG$-modules. Our aim is to illustrate that the predicted equalities follow naturally from highest weight  theory in the case under consideration.

\medskip
Let $\sigma: \mathbb{F} \to \mathbb{F}$ be the Frobenius automorphism of $\mathbb{F}$ with $f\sigma = f^{2}$ for all $f \in \mathbb{F}$. Let $S$ denote the natural module for $G$, realized over the field $\mathbb{F}_{2^{n}} \subseteq \mathbb{F}$. We note that the natural module $S$ is self-dual in this case. Let $I = \{0,1,2,\ldots, n-1\}$. By Steinberg's tensor product theorem, the
module $$S \otimes S^{\sigma} \otimes \ldots \otimes S^{\sigma^{n-1}}$$ is the Steinberg module ${\rm St}$, which is projective and simple.

\medskip
This allows us to describe all the simple $\mathbb{F}G$-modules up to isomorphism (there are $2^{n}$ of these, which is the number of $2$-regular conjugacy clases of $G$). For each non-empty subset $J$ of $I$, let 
$$S_{J} = \otimes_{j \in J} S^{\sigma^{j}},$$ and let $S_{\emptyset}$ denote the trivial module $\mathbb{F}.$ Then each $S_{J}$ is simple, since $S_{I}$ is   simple, and no two of these tensor products are isomorphic on consideration of their Brauer characters. For each subset $J$ of $I$, denote  $I \setminus J$ by $ J^{\prime}.$ We consider $S_{J}$ as the $m(J)$-th  simple module, where $m(J) = \sum_{j \in J} 2^{j}$ (and $m(\emptyset) =0$).

\medskip
We claim that the projective cover of $S_{\emptyset}$ is ${\rm Ind}_{T_{+}}^{G}(\mathbb{F}),$ where $T_{+}$ is a non-split torus of $G$ of order $2^{n}+1$, and that for each non-empty proper subset $J$ of $I$, the projective cover of $S_{J}$ is $S_{J^{\prime}} \otimes {\rm St}.$

\medskip
Recalling that each $S_{J}$ is self-dual, and noting that $$S_{J} \otimes S_{J^{\prime}} \otimes {\rm St} \cong {\rm St}\otimes {\rm St},$$ which has the projective cover of the trivial module as a summand, we conclude that the projective cover of $S_{J}$ is a summand of the projective module 
$S_{J^{\prime}} \otimes {\rm St}$ for each proper non-empty subset $J$ of $I$.

\medskip
Also the projective cover of the trivial module is a summand of the module $${\rm Ind}_{T_{+}}^{G}(\mathbb{F}),$$  which has dimension  $2^{2n}-2^{n}.$

\medskip 
It follows that we have $$|G| = \sum_{J \subseteq I} {\rm dim}_{\mathbb{F}}(S_{J}) {\rm dim}_{\mathbb{F}}(P(S_{J})) \leq (2^{2n}-2^{n}) + (2^{n}-1)2^{2n} = (2^{n}+1)2^{n}(2^{n}-1) = |G|.$$ Hence we must indeed have $$P(S_{J}) \cong S_{J^{\prime}} \otimes {\rm St}$$ for each non-empty subset $J$ of $I$, and $$P(\mathbb{F}) \cong {\rm Ind}_{T_{+}}^{G}(\mathbb{F}),$$ as claimed.

\medskip
In particular, we have $$S_{J} \otimes P(S_{J}) \cong {\rm St} \otimes {\rm St}$$ for each non-empty subset $J$ of $I$. We note also that $P(\mathbb{F})$ is a summand of ${\rm St} \otimes {\rm St},$ and by consideration of the dimensions of the projective indecomposable $\mathbb{F}G$-modules we must have $${\rm St} \otimes {\rm St} \cong P(\mathbb{F}) \oplus {\rm St}.$$

\medskip
Lifting to projective $RG$-modules, we may conclude that $$T_{c,2,G} = 2^{n} P(R) + (2^{n}-1) {\rm St}$$ in the Green ring (with obvious abuse of notation).

\medskip
Now let us recall the general formula $$T_{c,2,G} \cong \sum_{ x \in G_{2^{\prime}}/G} {\rm Ind}_{C_{G}(x)}^{G}(P_{1,2,C_{G}(x)}).$$

\medskip
If $x$ is a non-identity element of the non-split torus $T_{+}$ of order $2^{n}+1$, then we have $$P_{1,2,C_{G}(x)} \cong R$$ and $C_{G}(x) \cong T_{+}.$ There are $2^{n-1}$ conjugacy classes of such elements $x$, and we obtain $P_{x,2,G} \cong P(R)$ for each such $x$.

\medskip
If $T_{-}$ is a split torus of $G$ of order $2^{n}-1$ and $y$ is a non-identity element of $T_{-}$, then since there are $2^{n-1}-1$ conjugacy classes of such $y$, each with $P_{1,2,C_{G}(y)} \cong R,$ we must conclude that $$P_{y,2,G} \cong P(R) \oplus {\rm St} \oplus {\rm St}$$ for each such $y$.
This now forces $$P_{1,2,G} \cong P(R) \oplus {\rm St}.$$ We could also deduce this directly from the fact that $P(R)$ occurs with multiplicity one when $P_{1,2,G}$ is expressed a $\mathbb{Z}$-linear combination of projective indecomposables in the Green ring of $RG$. Note that each non-identity conjugacy class of elements of odd order in $G$ intersects  $T_{+} \cup T_{-}$  in a pair of mutually inverse elements.

\medskip
Notice that we now have $P_{y,2,G} - P_{1,2,G} = {\rm St}$ for each non-identity $y \in T_{-}$ and $P_{1,2,G} - P_{x,2,G} = {\rm St}$ for each non-identity $x \in T_{+}.$ Hence the Steinberg module ${\rm St}$ is detected by the modules $P_{z,2,G}$ for $2$-regular $z \in G$, and we have explicitly confirmed that each such $P_{z,2,G}$ really is a module (in this situation, the only case which is not immediately transparent is $P_{1,2,G}).$

\medskip
Let us examine the Lefschetz conjugation module in this context. The description of the simple $\mathbb{F}G$-modules tell us that the Steinberg module is the only projective simple $\mathbb{F}G$-module. Hence $RG$ has only one $2$-block of defect zero. The Sylow $2$-subgroup $U$ of $G$ is TI in this case, (by inspection, or by general theory). Since defect groups are always Sylow intersections, every block of positive defect of $RG$ has defect group $U$, and by a Theorem of Fong, [4], $RN_{G}(U)$ has a unique block (necessarily with defect group $U$), so by Brauer's First Main Theorem, the principal block is the only block of $RG$ with defect group $U.$ Hence the blocks of $RG$ are the principal block $B$ and the block $B_{s}$ containing the 
Steinberg module ${\rm St}.$

\medskip
Let us examine $T_{c,2,B}$ in the Green ring of $RG$ using our general formula  (here in the case $p =2$),
$${\rm T}_{c,p,B} = \sum_{Q/G} {\rm Ind}_{N_{G}(Q)}^{G}( P[\mathcal{L}_{c}(B(Q))]).$$ As noted in [15],  we may restrict the $p$-subgroups $Q$ on the right hand side to those $Q$ which satisfy $Q = O_{p}(N_{G}(Q))$, for if $O_{p}(N_{G}(Q)/Q) \neq 0,$ we have $\mathcal{L}_{c}(B(Q)) = 0.$ 
For our present group $G$ then,  we need only consider $Q =1 $ and $Q = U.$ We know from [8] that $\mathcal{L}_{c}(B)$ is already a virtual projective,  and as $N_{G}(U)/U$ is a cyclic group of order $2^{n}-1$, we see that $\mathcal{L}_{c}(B(U))$ is just the direct sum of $2^{n}-1$ copies of the trivial module $R$ (recalling that $B(U) = RN_{G}(U)/U$ in this case). The projective cover of the trivial module of $N_{G}(U)$ is equal to
${\rm Ind}_{T_{-}}^{N_{G}(U)}(R).$ 

\medskip
We saw in the earlier analysis that in the Green ring of $RG$, we have $${\rm Ind}_{T_{-}}^{G}(R) = 2{\rm St} + P(R).$$
Hence our formula simplifies in the Green ring to $$T_{c,2,B} = [B - {\rm Ind}_{N_{G}(U)}^{G}(RN_{G}(U))] + (2^{n}-1)[2{\rm St}+P(R)].$$

\medskip
We proved earlier that in the Green ring for $RG$, we have $$T_{c,2,G} = 2^{n}P(R) + (2^{n}-1){\rm St}$$ and $${\rm St} \otimes {\rm St} = P(R)+{\rm St}.$$ Since $$T_{c,2,G} - T_{c,2,B} = {\rm St} \otimes {\rm St}  = {\rm St} +P(R)$$ (recalling that ${\rm St}$ is self-dual), we thus obtain that 
$$T_{c,2,B} =( 2^{n}-1)P(R) + (2^{n}-2){\rm St},$$ and our module formula above now simplifies (after some rearrangement) in the Green ring to
 $$B  + 2^{n}{\rm St} = {\rm Ind}_{N_{G}(U}^{G}(RN_{G}(U))$$ (recall that $B$ and $RN_{G}(U)$ are modules under conjugation action for $G$ and $N_{G}(U)$ respectively). Hence in the Green ring for $RG$, we have $$\mathcal{L}_{c}(B) = -2^{n}{\rm St},$$ so that the unique block $B$ of positive defect of $RG$ is a block of Lefschetz  type. Also, $$\mathcal{L}_{c}(B_{s}) = P(R)+{\rm St},$$ yielding 
$$\mathcal{L}_{c}(RG) = P(R)-(2^{n}-1){\rm St}.$$ Then $\mathcal{L}_{c}(RG)$ is not the zero element of the Green ring of $RG,$ but does have ``virtual rank" zero. 

\medskip
\noindent {\bf Remark 3.1 :} While this example worked out very cleanly, in general, there are some limitations to the virtual characters $\Psi_{x,p,G}$ for the general finite group $G$ as $x$ runs over representatives for $p$-regular  conjugacy classes. For example, these generalized characters are all rational-valued, so if $G$ has any virtual projective characters  whose values are not rational, then $\{\Psi_{x,p,G} \}$ can not form a $\mathbb{Z}$-basis for the $\mathbb{Z}$-module of virtual characters of virtual projective $RG$-modules. In particular, these virtual characters will not distinguish between different algebraically conjugate characters in $p$-blocks of defect zero.

\medskip
For example, if we consider the case $G \cong A_{5}, p = 3,$ we find that $\Psi_{1,3,G}$ takes value $21$ on $1_{G}$ and value $1$ on all non-identity elements of order prime to $3$, while we find that for an involution $t$ of $G$, $\Psi_{t,3,G}$ takes value $15$ on $1_{G}$, value $3$ on conjugates of $t$, and value $0$ on elements of order $3$. It then follows that $\Psi_{1,3,G} - \Psi_{t,3,G}$ is a generalized character of norm squared $2$, and in fact we have 
$$\Psi_{1,3,G} -  \Psi_{t,3,G} = \chi + \chi^{\sigma},$$ where $\chi$ is an irreducible character of degree $3$ and $\chi^{\sigma}$ is its other algebraic conjugate. These are the only (irreducible characters in) blocks of defect zero of $RG$. In this case, it might be argued that while the $\Psi_{z,3,G}$ are characters, but do not separate the two algebraically conjugate $3$-blocks of defect zero, they do detect the fact that there are two $3$-blocks of defect zero. Here, we knew already from the general results of [16] that the generalized characters $\Psi_{z,3,G}$ are indeed characters, but it can also be seen directly in this small case.

\medskip
Another example where this approach bears fruit is when $p = 5$ and $$G \cong A_{5} \cong {\rm PSL}(2,5).$$ We omit the details, but letting $t$ be an involution  of $G$ and $x$ be an element of order $3$, we find that $\Psi_{1,5,G} - \Psi_{x,5,G}$ is the Steinberg character of degree $5$.

\medskip
We examine the similar case with $G \cong {\rm PSL}(2,7)$ and $p = 7$ in more (but not exhaustive) detail. In this case,  a set of representatives  for the conjugacy classes of $7$-regular elements is $\{1,t,u,x \}$ where $t$ has order $2$, $u$ has order $4$, and $x$ has order $3$. There are four projective indecomposable $RG$ modules: these are $P_{1} \cong {\rm Ind}_{S_{4}}^{G}(R)$ ( it does not matter which conjugacy class of $S_{4}$'s we choose),
${\rm St}$, which is the non-trivial (projective) summand of ${\rm Ind}_{F_{21}}^{G}(R)$ (and $F_{21}$ denotes the Borel subgroup, which is a Frobenius group of order $21$), $P_{3}$, which is a projective indecomposable $RG$-module of rank $14$ whose character is the sum of the irreducible character of $G$ of degree $8$ and the irreducible character of $G$ of degree $6$, and $P_{4}$, whose character is the sum of the irreducible character of degree $8$ of $G$ and the two exceptional irreducible characters of degree $3$ in the principal $7$-block of $G$.

\medskip
We deduce from the character table of $G$ that $$P_{1,7,G} \cong {\rm St} \otimes {\rm St} \cong P_{1} \oplus {\rm St} \oplus {\rm St} \oplus P_{3}\oplus P_{4},$$ that $$P_{t,7,G} \cong {\rm P}_{1} \oplus P_{3},$$ that $$P_{x,7,G} \cong P_{1} \oplus {\rm St} \oplus {\rm St} \oplus {\rm St} \oplus  P_{3}\oplus P_{4},$$ and that $$P_{u,7,G} \cong P_{1} \oplus {\rm St} \oplus P_{3} \oplus P_{4}.$$

\medskip
Hence each $\Psi_{z,7,G}$ is the character of a projective $RG$-module, and we have that (the  character of) 
$${\rm St}  = \Psi_{1,7,G} -\Psi_{u,7,G} = \Psi_{x,7,G} - \Psi_{1,7,G}.$$ It now follows that the $\mathbb{Z}$-submodule of the Green ring generated by 
$$\{P_{1,7,G},P_{t,7,G},P_{x,7,G}, P_{u,7,G} \}$$ has $\mathbb{Z}$-basis $\{{\rm St},P_{4},P_{1} + P_{3} \}.$

\medskip
Let us also examine the Lefschetz conjugation modules for the principal block $B$ of $G$ in the cases $G \cong A_{5}$ and $G \cong {\rm PSL}(2,7)$ (both for $p =3$, to avoid Lie-type characteristic $p$ cases). In the case $G \cong A_{5},$ we see that the virtual character $\mathcal{L}_{c}(B)$, as well as vanishing on elements of order $3$, takes value $2$ on elements of order $5$, value $-18$ on $1_{G},$ and value $-2$ on involutions. Hence the trivial character occurs with multiplicity zero, the irreducible character of degree $5$ occurs with multiplicity $-2$, the irreducible character of degree $4$ occurs with multiplicity $-2$, and each irreducible character of degree $3$ occurs with multiplicity zero. Since virtual projectives in the Green ring are determined by their virtual characters, it follows that $$\mathcal{L}_{c}(B) = 2{\rm St}_{3}$$ in the Green ring, since for $U \in {\rm Syl}_{3}(G)$ we see that  $${\rm St}_{3} = R - {\rm Ind}_{N_{G}(U)}^{G}(R)$$ has virtual character $-(\alpha + \beta),$ where $\alpha$ is irreducible of degree $4$ and $\beta$ is irreducible of degree $5,$  so that ${\rm St}_{3}$ is the negative of a projective $RG$-module.

\medskip
In the case $G \cong {\rm PSL}(2,7),$ (and $p=3$), the irreducible characters in the principal block $B$ of $RG$ have degree $1,7,$ and $8$. No element of order $4$ or $7$ lies in $N_{G}(U)$ for $U \in {\rm Syl}_{3}(G)$. Hence the virtual character $\theta$ afforded by $\mathcal{L}_{c}(B)$ takes value $-54$ on $1_{G},$ value $2$ on elements of order $4$, value $2$ on elements of order $7$, value $0$ on elements of order $3$, and 
value $-6$ on involutions. It follows that $\theta$ contains the trivial character with multiplicity $0$, the degree $6$ irreducible with multiplicity $-4$, the degree $8$ irreducible with multiplicity $-2$, the  degree $7$ irreducible with multiplicity $-2$ and each irreducible of degree $3$ with multiplicity zero.

\medskip
 Since all non-principal blocks of $RG$ have defect zero, and the two projective indecomposable $B$-modules have respective ranks $9$ and $15$, we see that $$\mathcal{L}_{c}(B) = -4P_{6} -2P_{15}$$ in the Green ring of $RG$, where $P_{6}$ is the unique projective indecomposable of rank $6$ of $RG$ and $P_{15}$ is the unique projective indecomposable $RG$-module of rank $15$.

\medskip
The example $G = {\rm GL}(3,2), p = 2$ is also somewhat illuminating. We consider it because it is an example for which the Sylow $p$-normalizer is not strongly $p$-embedded, unlike the cases considered up to now. There are four projective indecomposable $RG$-modules, which are the projective cover of the trivial module $P(R)$, of rank $8$, which is induced form the trivial module of the Sylow $7$-normalizer, the Steinberg module of rank $8$, which may be induced from a non-trivial module of rank one of the Sylow $7$-normalizer, and two projective indecomposables of rank $16$, which are the projective covers of the natural module for $G$, and the projective cover of the dual of the natural module, labelled $P_{n}$ and $P_{d}$ respectively. If $\alpha$ is the character of the (lift of) the natural module, then the character afforded by $P_{n}$ is $\alpha + \gamma+ \chi$, where $\gamma$ is the unique irreducible character of degree $6$ and $\chi$ is irreducible of degree $7$. The character afforded by $P_{d}$ is the complex conjugate of the character of $P_{n}.$

\medskip
Let $\theta $ denote the virtual projective character afforded by $\mathcal{L}_{c}(B)$, where $B$ is the principal block of $RG$. Since $\theta$ vanishes on all non-identity elements of even order, we only need to concern ourselves with $\theta(1),\theta(x), \theta(y)$ and $\theta(y^{-1})$, where $x$ has order $3$ and $y$ has order $7$. Note that $$\theta(1) = 104 - 168-168 + 168 =-64.$$

\medskip
Since $y$ normalizes no non-trivial $2$-subgroup of $G$, we see that $$\theta(y) = \theta(y^{-1}) = 6.$$ Since $x$ lies in two of the subgroups of $G$ isomorphic to $S_{4}$ and $|C_{H}(x)| = 3$ for each such subgroup  $H$, we find that $\theta(x) = -4$. We check easily from this information and the character table for $G$ that $$\theta = -2( \alpha + \overline{\alpha} + 2\gamma + 2 \chi),$$ which is the virtual character afforded by the virtual projective 
 $-2(P_{n}+P_{d}).$ Since virtual projectives are determined in the Green ring by their virtual character, we have $$\mathcal{L}_{c}(B) = -2(P_{n}+P_{d})$$ in the Green ring of $RG$.

\medskip
\section{The Frobenius automorphism of $\mathbb{F}$ and the virtual characters $\Psi_{1,p,G}$}

\medskip
In general, it is not necessarily the case that the Galois group $\mathbb{Q}[\omega]/\mathbb{Q}$ acts on the set of Brauer irreducible characters associated to the chosen $p$-modular system $(R,\mathbb{K},\mathbb{F})$ when $\omega$ is a complex primitive $|G|$-th root of unity. However, if $\mathbb{F}_{p^{n}}$ of size $p^{n}$  is a subfield of $\mathbb{F}$ which is a splitting field for $G$, then the cyclic Galois group ${\rm Gal}(\mathbb{F}_{p^{n}}/\mathbb{F}_{p})$ of order $n$ does permute the irreducible Brauer characters associated to this $p$-modular system, and also acts in a compatible fashion on the characters of the projective indecomposable $RG$-module. In fact (as proved in K\"ulshammer [9], for example), the smallest possible choice of $p^{n}$ is the smallest power of $p$ such that $y^{p^{n}}$ is $G$-conjugate to $y$ for each $p$-regular $y \in G$. 

\medskip
To be precise, if we let $\sigma$ denote the Frobenius automorphism of $\mathbb{F}$, then $\sigma$ acts on the set of Brauer irreducible characters $\{ \phi_{i} : 1 \leq i \leq \ell(G) \}$,  via $\phi_{i}^{\sigma}= \phi_{j},$ where $\phi_{j}(y) = \phi_{i}(y^{p})$ for each $p$-regular $y \in G$. Also $\sigma$ acts on the set $\{ \theta_{i} \}$ of characters of projective indecomposable $RG$ modules in a similar manner, and we have 
$$\langle \theta_{i}^{\sigma}, \phi_{j}^{\sigma} \rangle = \langle \theta_{i}, \phi_{j} \rangle = \delta_{ij}$$ for all $i,j,$ where $\langle, \rangle$ denotes the usual character inner product.

\medskip
In fact, for each subgroup $H$ of $G$, we see that the Frobenius automorphism $\sigma$ of $\mathbb{F}$ acts on the Brauer irreducible characters and the characters of the projective indecomposable $RH$-modules of each subgroup $H$ of $G$ (even though $H$ might itself require a larger splitting field than $G$ does). This action is compatible in the sense that whenever $\phi$ is a $\sigma$-stable virtual Brauer character of $G$, then all Brauer irreducible characters of $H$ in the same $\sigma$-orbit occur with equal multiplicity in ${\rm Res}_{H}^{G}(\phi).$

\medskip
When the rational-valued virtual projective characters $\Psi_{x,p,G}$ are decomposed as $\mathbb{Z}$-linear combinations of $\{ \theta_{i} \}$, it is the case that for each $i$, all elements of the $\langle \sigma \rangle$-orbit of $\theta_{i}$ occur with the same coefficient as $\theta_{i}$ does, and this coefficient is $\langle \Psi_{x,p,G}, \phi_{i} \rangle$ (recall that this inner product makes sense, as although $\phi_{i}$ is only defined on $p$-regular elements, we know that $\Psi_{x,p,G}$ vanishes on all $p$-singular elements).

\medskip
It may happen that $\Psi_{1,p,G}$ is induced from a virtual projective character of a normal subgroup $N$ of $G$, but it is not necessarily the case when that happens that the virtual projective of $N$ is itself rational-valued (in particular, it may be different from $\Psi_{1,p,N}$). This may sometimes assist in proving that $P_{1,p,G}$ is a genuine projective module. We illustrate with an example.

\medskip
Let $G = S_{5}$ and take $p=2.$ If $P_{1,2,G}$ is a projective $RG$-module, then by Green's indecomposability theorem, it must be induced from a projective $RN$-module $M$, where $N = A_{5}.$ Since $P_{1},$ the projective cover of the trivial module $R$ of $G$ must occur with multiplicity one as a summand of $P_{1,2,G},$ we see that the projective cover of the trivial $RN$-module must occur with multiplicity one as a summand of $M$. Again, using Green's indecomposabilty theorem (and the fact that the projective cover of the trivial $RN$-module has rank $12$) shows that the projective cover of the trivial $RG$-module has $R$-rank $24$. Now the virtual projective generalized character $\Psi_{1,2,G}$ takes value $56$ at the identity, value $1$ on elements of order $5$, and value $2$ on elements of order $3$ (as well as vanishing on all $2$-singular elements).

\medskip
It follows that if $P_{1,2,G}$ is a projective $RG$-module, then the $RN$-module $M$ it is induced from can not be both of $G$-stable isomorphism type and have rational-valued character (otherwise, all of its values would be even integers). We note (using for example, the fact that $N \cong {\rm SL}(2,4)$) that there are only two projective indecomposable  $RN$-modules whose characters are not rational-valued, and that these are interchanged by the Frobenius automorphism of $\mathbb{F}$ by consideration of their values on elements of order $5$). These two modules each have $R$-rank $8$, and their direct sum has rational-valued character, and is of $G$-stable isomorphism type.

\medskip 
It now follows that $M$ (which must have $R$-rank $28$) can not be the direct sum of the projective cover of the trivial-module with the direct sum of these two rank $8$ projective $RN$-modules (for in that case, $M$ would be of $G$-stable isomorphism type, and have rational-valued character). The only possibilities are that $M$ is the direct sum of the projective cover of the trivial module with two copies of the rank $4$ projective indecomposable $RN$-module and one of the rank $8$ projective indecomposable modules. Each of these choices induce the same way to $G$, and the result affords character 
$\Psi_{1,2,G}.$

\medskip
Hence $P_{1,2,G}$ is indeed a projective $RG$-module in this case, and the point of the example is to illustrate that this forces the existence of two non-isomorphic projective indecomposable $RN$-modules which are interchanged by the Frobenius automorphism of $\mathbb{F}$.

\medskip
To be explicit about $P_{1,2,G}$ we may conclude that it is isomorphic to the direct sum of the unique projective cover of the trivial module (which has rank $24$) with two copies of the rank $8$ projective indecomposable $RG$-module  (induced from the unique rank $4$ projective indecomposable of $N$) and the unique rank $16$ projective indecomposable $RG$-module. This differs from ${\rm Ind}_{N}^{G}(P_{1,2,N})$ by the addition of the rank $16$ projective indecomposable summand and one copy of the rank $8$ projective indecomposable summand. This difference is necessary to ``correct" the fact that ${\rm Ind}_{N}^{G}(\Psi_{1,2,N})$ only has degree $32$ and takes value $2$ on all elements of order $5$.

\section{ Remarks on the case $G = S_{n}, p = 2$}

\medskip
We have noted that when $B$ is a block of defect zero of $RG$, then $\mathcal{L}_{c}(B)$ is a projective $RG$-module which contains the projective cover of the trivial module as a summand of multiplicity one. In all the cases studied so far, it has been the case that when $B$ is a block of positive defect of $RG$, then the virtual module $\mathcal{L}_{c}(B)$ has the form $-M$ in the Green ring, where $M$ is a genuine projective $RG$ module, and $M$ does not have the projective cover of the trivial module as a summand. We will show by example that this is not always the case.

\medskip
If $G$ is a finite group  such that $RG$ has no block of defect zero, (so that $G$ must definitely have order divisible by $p$), we have seen that $\mathcal{L}_{c}(RG)$ has virtual rank zero, while all blocks $B$ of $RG$ have positive defect, so there is at least one block $B$ of $RG$ such that $\mathcal{L}_{c}(B)$ has non-negative virtual rank. Hence we now seek a group $G$ such that $\mathcal{L}_{c}(RG)$ is not zero in the Green ring and $RG$ has no block of defect zero.

\medskip
We note an easy general fact (this applies a weaker version  of a result in [14] which was proved in the context of Quillen's conjecture on the $p$-subgroup complex, but is applied here in the context of attaching a non-trivial coefficent system to the $p$-subgroup complex):

\medskip
\noindent {\bf Lemma 5.1 :} \emph{Suppose that the finite group $G$ contains a (necessarily $p$-regular) element $y$ which normalizes no non-trivial $p$-subgroup of $G$. Then the character afforded by $\mathcal{L}_{c}(RG)$ takes value $|C_{G}(y)|$  at $y$. In particular, there is at least one block $B$ of $RG$ such that the  value of the character afforded by $\mathcal{L}_{c}(B)$ takes a value at $y$ which lies outside $J(R).$ Hence $\mathcal{L}_{c}(RG)$ and $\mathcal{L}_{c}(B)$ are (possibly equal) non-zero elements of the Green ring of $RG$.}

\medskip
\noindent {\bf Proof:} Notice that the hypotheses of the Lemma imply that $O_{p}(G) = 1,$ and that $y$ centralizes no non-identity $p$-element of $G$, so that $|C_{G}(y)|$ is prime to $p$.  Since no conjugate of $y$ lies in $G_{\sigma}$ for any non-empty chain $\sigma \in \mathcal{S}_{p}(G)$, we see that only the character afforded by the conjugation action on $RG$ contributes to the value at $y$ of the character afforded by $\mathcal{L}_{c}(RG).$
That contribution is $|C_{G}(y)|$ which is coprime to $p$ and hence certainly non-zero. Hence $\mathcal{L}_{c}(RG)$ is certainly not zero in the  Green ring of $RG$. Clearly there must be a block $B$ of $RG$ such that the value taken at $y$ by the character afforded by  $\mathcal{L}_{c}(B)$ lies outside $J(R),$ and is hence also non-zero.

\medskip
\noindent {\bf Corollary 5.2:} \emph{ If $RG$ has no block of defect zero, and $G$ contains a non-identity element $y$ which normalizes no non-identity $p$-subgroup of $G$, then $RG$ has least one block $B$ such that $\mathcal{L}_{c}(B)$ is not a genuine projective module, and at least one block $B^{\prime}$ (possibly equal to $B$), such that $\mathcal{L}_{c}(B^{\prime})$ is not  the negative of a genuine projective module. In particular, if $RG$ has only one block, then $\mathcal{L}_{c}(RG)$ is (in the Green ring of $RG$), neither a projective module, nor the negative of a projective module.}

\medskip
\noindent {\bf Remark 5.3:} The group $G = S_{5}$ has no $2$-block of defect zero, while a $5$-cycle of $S_{5}$ normalizes no non-identity $2$-subgroup of $G$, so the first part of the Corollary applies in the case $G = S_{5}$ and $p = 2.$ 

\medskip
We will now examine some more general statements about the Lefschetz conjugation modules of $2$-blocks of $S_{n}.$ The next Lemma is elementary, but we include a proof for the convenience of the reader.

\medskip
\noindent {\bf Lemma 5.4:} \emph{For every odd integer $n$, an $n$-cycle of $S_{n}$ normalizes no non-trivial $2$-subgroup of $S_{n}$. If $n$ is an even integer, then an $(n-1)$-cycle of $S_{n}$ normalizes a non-trivial $2$-subgroup of $S_{n}$ if and only if $n$ is a power of $2$.}

\medskip
\noindent {\bf Proof:} Suppose first that $n$ is odd. Let $\sigma = (12 \ldots n-1 n) \in S_{n}$ without loss of generality, and let $T = \langle \sigma \rangle$. Suppose that $A$ is a non-trivial $2$-subgroup of $S_{n}$ of minimal order subject to being normalized by $T$. Then $A$ has an orbit of odd length, so $A$ must fix a point. Then $T$ permutes the fixed-points of $A$, and $T$ is transitive, so $A$ must fix all points, a contradiction.

\medskip
Suppose now that $n$ is even. Let $\sigma = (12 \ldots n-1) \in S_{n}$ without loss of generality, and let $H = \langle \sigma \rangle$. Suppose that $A$ is a non-trivial $2$-subgroup of $S_{n}$ of minimal order subject to being normalized by $H$. Then $A$ is an elementary Abelian $2$-group. Let $X = AH$. Then $X$ is clearly transitive, for $A$ must move the point $n$ (otherwise $AH$ is isomorphic to a subgroup of $S_{n-1}$ and $n-1$ is now odd, so the argument above applies). Then $$AH = HA= {\rm Stab}_{AH}(n)A,$$ so that $A$ is itself transitive. Hence $n$ divides $|A|$, and $n$ must be a power of $2$.

\medskip
On the other hand, if $n = 2^{a}$ is a power of $2,$ then $ F = \mathbb{F}_{2^{a}}$ is (additively) an elementary Abelian group of order $2^{a}$, and its  (cyclic) multiplicative group $C$ of order $2^{a}-1$ acts regularly on the non-zero elements of $F$ by multiplication. Then $FC$ is a doubly transitive subgroup of ${\rm Sym}(F) \cong S_{n}$ with $C$, the (cyclic) stabilizer of $0$, acting as $n-1$ cycle on the remaining elements of $F$.

\medskip
\noindent {\bf Theorem 5.5:} \emph{ Let $n$ be an integer which is neither a power of $2$, nor a triangular number, and set $G = S_{n}, p =2$.
Then the (virtually projective) Lefschetz conjugation module $\mathcal{L}_{c}(RG)$ is non-zero, and is neither a module, nor the negative of a module, in the Green ring of $RG$.  Furthermore, $RG$ has a $2$-block $B$ such that the virtual projective module $\mathcal{L}_{c}(B)$  is not a module in the Green ring of $RS_{n}$, and a block $B^{\prime}$ (possibly equal to $B$), such that $\mathcal{L}_{c}(B^{\prime})$ is not the negative of a module in the Green ring of $RG$ .}

\medskip
\noindent {\bf Remark 5.6:} For example, when $p=2$ and $G = S_{5},$ we find that the Lefschetz conjugation module of $RG$ has character which takes value $0$ at $1_{G}$, takes value $5$ at all $5$-cycles, and value $-3$ at all elements of order $3$ (and vanishes on all elements of even order). It is fairly straightforward to check from this information that  in the Green ring of $RG$, we have $\mathcal{L}_{c}(B) = 3P_{16}$  (a genuine module) when $B$ is the principal $2$-block, where $P_{16}$ is the unique projective indecomposable $RG$ module of rank $16$, and we have $L_{c}(B^{\prime}) = -2P_{s} -2P_{16},$ (the negative of a module) when $B^{\prime}$ is the $2$-block of $RG$ of defect one, and $P_{s}$ is the unique projective indecomposable module of $RG$ of rank $8$. Note the fact that, as expected, the projective cover of the trivial $RG$-module $P(R)$ appears with zero multiplicity in each of these virtual modules.

\medskip
For the sake of completeness, we would like to exhibit an example of a group $G$ and a single block $B$ of $RG$ such that in the Green ring of $RG$, the Lefschetz conjugation module $\mathcal{L}_{c}(B)$ is a non-zero virtual module of ``mixed type" (neither a module nor the negative of a module). There may be smaller examples, but one example is provided by:

\medskip
\noindent {\bf Remark 5.7:} It is well-known that the sporadic finite simple groups $M_{22}$ and $M_{24}$ each have a unique $2$-block (see, for example, (Harris, [7])). Hence $M_{24}$ with $p=2$ is an example of a finite group $G$ which satisfies $O_{p}(G) = 1$ and for which $RG$ only has one $p$-block.

\medskip
In this case, we may note from the list of maximal subgroups of $M_{24}$ in the Atlas, [3], that an element of order $23$ in $G$ normalizes no non-trivial $2$-subgroup of $G$. To show this, the degree $23$ permutation representation of the maximal subgoup $M_{23}$  can be used in the same fashion as in 5.4 above. Since $RG$ is the principal $2$-block of $G$ in this case, and (using 5.5), $\mathcal{L}_{c}(RG)$ is a non-zero virtual module of virtual rank $0$, its Lefschetz conjugation module is a non-zero element in the Green ring which is neither a module nor the negative of a module.

\section{Some relevant $R[G \times G]$-modules}

\medskip
One motivation for this section is that when we consider $RG$ as an $R[G \times G]$-module ( cf Green [5]) in now standard fashion, the indecomposable summands are precisely the blocks of $RG$, and the projective summands are precisely the blocks of defect zero of $RG$. Hence it might be profitable to consider $RG$ in the stable module category of $R[G \times G].$

\medskip
We begin with a simple, yet suggestive, observation:

\medskip
\noindent {\bf Theorem 6.1:} \emph{ The virtual $R[G \times G]$-module $$\sum_{ \sigma \in  \mathcal{S}_{p}(G)/G} (-1)^{|\sigma|} 
{\rm Ind}_{G_{\sigma} \times G_{\sigma}}^{G \times G}(RG_{\sigma})$$ is a virtual projective in the Green ring for $R[G \times G]$. In particular, this virtual module is zero if the p-subgroup complex $\mathcal{S}_{p}(G)$ is contractible.}

\medskip
\noindent {\bf Proof:} We first identify $G$ with the diagonal subgroup $\Delta(G)$ of $G \times G$. We also recall that J.A. Green proved in [5] that as $R[G \times G]$-module, we have $$RG \cong {\rm Ind}_{\Delta(G)}^{G \times G}(R).$$ Noting that P.J. Webb has proved that
$$\sum_{\sigma \in \mathcal{S}_{p}(G)/G}(-1)^{|\sigma|} {\rm Ind}_{\Delta(G_{\sigma})}^{\Delta(G)}$$ is virtual projective as $R\Delta(G)$-module (and is zero if the $p$-subgroup complex $\mathcal{S}_{p}(G)$ is contractible), we may induce this virtual module to a virtual projective in the Green ring of $R[G \times G]$. Inducing the individual (virtual) summands first from $\Delta(G_{\sigma})$ to $G_{\sigma} \times G_{\sigma}$, then to $G \times G$, and using transitivity of induction, we obtain the Theorem as stated.

\medskip
A similar argument, starting from the results of [8], yields the following analogous statement when $B$ is a block of $RG$, viewed as an $R[G \times G]$--module.

\medskip
\noindent {\bf Theorem 6.2:} \emph{ Let $B$ be a block of $RG$. Then the virtual $R[G \times G]$-module $$\sum_{ \sigma \in  \mathcal{S}_{p}(G)/G} (-1)^{|\sigma|} 
{\rm Ind}_{G_{\sigma} \times G_{\sigma}}^{G \times G}(B_{\sigma} \otimes RG_{\sigma})$$ is a virtual projective in the Green ring for $R[G \times G]$. }

\medskip
The key observation in proving the latter statement is that when $M$ is an $R\Delta(G)$-module which extends to an $R[G \times G]$-module, then $${\rm Ind}_{\Delta(G)}^{G \times G}(M) \cong M \otimes RG$$ as $R[G \times G]$-module.  We apply this by noting that the conjugation action of $G$ on $B$ may be viewed as the restriction from $G \times G$ to $\Delta(G)$ of the usual $G \times G$-action.

\medskip
\noindent {\bf Remark 6.3:} Before progressing to the next section, we recall that the character of $G \times G$ afforded by its action $ a \mapsto x^{-1}ay$ for all $ a \in RG$ and $(x,y) \in G \times G$ is $$\sum_{ \chi \in {\rm Irr}(G)} \overline{\chi} \otimes \chi.$$

\section{The number of blocks of defect zero revisited}

\medskip
In the paper [12], we obtained a formula for the number of blocks of $RG$ with a given defect group $D$ of the finite group $G$ when $D$ is a normal $p$-subgroup of $G,$ in terms of group-theoretic information (by Brauer's First Main Theorem, the assumption that $D \lhd G$ is not really any loss of generality).
In any case, in the present context, we are primarily concerned with the case $D =\{1_{G}\},$ which is most relevant to Alperin's weight conjecture. 

\medskip
Let $\ell_{0} = \ell_{0}(G)$ denote the number of $G$-conjugacy classes of (necessarily $p$-regular) elements $y \in G$ such that $|C_{G}(y)|$ is not divisible by $p$. We note that the $\mathbb{F}$-span of the class sums of such elements span an ideal of $Z(\mathbb{F}G),$ which is customarily denoted by ${\rm Tr}_{1}^{G}(\mathbb{F}G).$ The ideal ${\rm Tr}_{1}^{G}(\mathbb{F}G)$ has a direct sum decomposition as a direct sum $J \oplus Z,$ where $J$ is the ideal consisting of nilpotent elements, and $Z$ is the $\mathbb{F}$-span of the idempotents in ${\rm Tr}_{1}^{G}(\mathbb{F}G).$ In fact, both $J$ and $Z$ are ideals of $Z(\mathbb{F}G),$ since ${\rm Tr}_{1}^{G}(\mathbb{F}G)$ annihilates $J(Z(\mathbb{F}G))$ (and even annihilates $J(\mathbb{F}G)$).

\medskip
For the case of counting blocks of defect zero of $RG$, it was shown in [12] that the number of blocks of defect zero is the rank of the matrix $[a_{ij}] : 1 \leq i, j \leq \ell_{0},$ where, for a Sylow $p$-subgroup $S$ of $G$, representatives $y_{i}$ and $y_{j}$ of classes of $p$-defect zero, we let $a_{ij}$ be the residue (${\rm mod}$  $p$) of the integer $$\frac{|\Omega_{ij}|}{|S|},$$ where $$\Omega_{ij} = \{ (z_{i},z_{j}): z_{i}^{-1} z_{j} \in S, z_{i} \in y_{i}^{G}, z_{j} \in y_{j}^{G} \}.$$

\medskip
By the Burnside-Frobenius class algebra constant formula (which uses the irreducible characters of $G$), this leads to 
$$\frac{|\Omega_{ij}|}{|S|} = \frac{1}{|C_{G}(y_{i})||C_{G}(y_{j})|}\left( \sum_{\chi \in {\rm Irr}(G)}\left(\frac{|G|}{\chi(1)}\right)\chi(y_{i}^{-1})\chi(y_{j})\langle {\rm Res}^{G}_{S}(\chi),1 \rangle\right).$$ We are only interested in $a_{ij},$ which is the residue $({\rm mod}$  $p)$ of this integer.

\medskip
Since $|C_{G}(y_{i})||C_{G}(y_{j})|$ is not divisible by $p$, and since $\frac{|G|}{\chi(1)}$ is always an integer, while $$\chi(y_{i}^{-1})\chi(y_{j})\langle {\rm Res}^{G}_{S}(\chi),1 \rangle$$ is always an algebraic integer, we may ignore terms for which $\frac{|G|}{\chi(1)}$ is divisible by $p$, that is, we may ignore those irreducible characters $\chi$ which do not have degree divisible by $|S|$. Hence the  only irreducible characters which we need to consider are those irreducible characters $\chi$ which lie in blocks of defect zero of $RG$. For such an irreducible character $\chi,$ we have 
$$\langle {\rm Res}^{G}_{S}(\chi),1 \rangle = \frac{\chi(1)}{|S|},$$ since $\chi$ vanishes on all $p$-singular elements of $G$.

\medskip
If there are no blocks of defect zero of $RG$, then we see directly that we have $a_{ij} = 0$ (in $\mathbb{Z}/p\mathbb{Z}$) for all $i,j$).
If there are any blocks of defect zero of $RG,$ then, letting $\Delta_{0}$ denote the sum of the blocks of defect zero of $RG,$  we have (after a little rearrangement), 
$$a_{ij} \equiv  \frac{[G:S]}{|C_{G}(y_{i})||C_{G}(y_{j})|} \sum_{ \chi \in {\rm Irr}(\Delta_{0})} \chi(y_{i}^{-1})\chi(y_{j}) ({\rm mod} \, p).$$

\medskip
\noindent {\bf Remark 7.1:} This last equation reflects the fact that the matrix $[a_{ij}]$ detects (mod $p$) (via group-theoretic information) the ``missing parts" of column inner products of the character table of $G$ coming from columns of $p$-regular elements whose centralizers have order prime to $p$, that is to say, those parts not detected by irreducible characters in blocks of $RG$ of positive defect.

\medskip
\noindent {\bf Remark 7.2:} We note that the matrix $E_{0}$, whose $(i,j)$-entry is $$\frac{|C_{G}(y_{j})|}{[G:S]} a_{ij}$$ is  the zero matrix if $RG$ has no block of defect zero, and is an idempotent matrix whose rank is the number of blocks of defect zero of $RG$ if $RG$ has a block of defect zero (the key point being that the matrix $[a_{ij}]$ is always defined, in purely group-theoretic terms, independently of any a priori knowledge of the existence (or otherwise) of blocks of defect zero of $RG$). The idempotent matrix $E_{0}$ represents the projection from the ideal ${\rm Tr}_{1}^{G}(\mathbb{F}G)$ onto its ``semisimple part" $Z$, with kernel $J$.

\medskip
It is a standard argument to verify that $$E = \lim_{n \to \infty}E_{1}^{p^{n}}$$ is an idempotent matrix (or the zero matrix) in $R$ whose trace is the number of blocks of defect zero of $RG$ for any preimage $E_{1}$ in $M_{\ell_{0} \times \ell_{0}}(R)$ of $E_{0},$ and that $E$ is also a preimage of $E_{0}$. 

\medskip
We remark, incidentally, that the matrix $[a_{ij}]$ (with entries in $\mathbb{F}$) is already idempotent (or zero) when $p = 2.$

\medskip
\noindent {\bf Remark 7.3:} For the purposes of illustration, we consider the case $G = {\rm PSL}(2,7)$ and $p =3$. There are four conjugacy classes of $3$-regular elements whose centralizers have order coprime to $3$ in $G$: these have representatives $y_{1}$ of order $2$, $y_{2}$ of order $4$, $y_{3}$ of order $7$   and $y_{4} = y_{3}^{-1},$ also of order $7.$ We have $|C_{G}(y_{i})| \equiv 1$ (mod $3$) for $i >1,$ and $|C_{G}(y_{1})| \equiv -1$ (mod $3$).

\medskip
In this case, the matrix $E_{0}$ is $$\left(\begin{array}{clcrc} 0 & -1 & 1 & 1\\ 1 & -1 & -1 &-1\\-1 &-1&-1&1\\-1&-1&1&-1 \end{array} \right),$$  which is  idempotent of rank $3$, reflecting the fact that $RG$ has three blocks of defect zero.

\medskip
Let us examine the matrix $E_{0}$ more carefully for a while in the general case. Let us label a set of representatives for the $p$-regular conjugacy classes of $G$ whose centralzers have order prime to $p$ by $y_{1},y_{2},\ldots, y_{\ell_{0}}$, and let us label the $f_{0}(G) = f_{0}$ irreducible characters of $G$ in blocks of defect zero of $G$ by $\chi_{1}, \chi_{2}, \ldots, \chi_{f_{0}}.$ Then $E_{0}$ is the image (mod $J(R)$) of the $\ell_{0} \times \ell_{0}$-matrix $E_{1}$ whose $(i,j)$-entry is $$\frac{1}{|C_{G}(y_{i})|}\left(\sum_{k=1}^{f_{0}} \chi_{k}(y_{i}^{-1})\chi_{k}(y_{j})\right).$$ We know that this matrix has maximum possible rank, which is $f_{0}.$ Notice that all entries of $E_{1}$ lie in $R$, and that $E_{1}$ has rational entries.
 
\medskip
Let $L$ be the $f_{0} \times \ell_{0}$-matrix whose $(k,j)$-entry is $\chi_{k}(y_{j})$, and let $M$  be the $\ell_{0} \times f_{0}$ matrix whose $(i,k)$-entry is $$\frac{1}{|C_{G}(y_{i})|}\chi_{k}(y_{i}^{-1}).$$ Then we have $ML = E_{1}.$ Let $L^{\prime}$ and $M^{\prime}$ denote the respective images (mod $J(R)$) of $L$ and $M$, so that $M^{\prime}L^{\prime}$ is the idempotent $\ell_{0} \times \ell_{0}$ matrix $E_{0}$ of rank $f_{0}.$

\medskip
We claim that we in fact have $L^{\prime}M^{\prime} = I_{f_{0} \times f_{0}}.$

\medskip
The $(i,j)$- entry of $LM$ is  $$\sum_{k=1}^{\ell_{0}} \frac{\chi_{i}(y_{k}) \chi_{j}(y_{k}^{-1})}{|C_{G}(y_{k})|}.$$ We claim that this is congruent to 
$\delta_{i,j}$ (mod $J(R)$), which will establish the claim.

\medskip
Since $\chi_{i}$ and $\chi_{j}$ both lie in blocks of defect zero of $RG$, we know that 
$$\sum_{y \in G_{p^{\prime}}/G}\frac{\chi_{i}(y)\chi_{j}(y^{-1})}{|C_{G}(y)|} = \delta_{i,j}.$$

\medskip
However, when $y$ is $p$-regular with $|C_{G}(y)|$ divisible by $p$, we have $$[G:C_{G}(y)]\frac{\chi_{i}(y)}{\chi_{i}(1)} \in R$$ and 
$$[G:C_{G}(y)]\frac{\chi_{j}(y^{-1})}{\chi_{j}(1)} \in R. $$ Since $[G:C_{G}(y]$ is not divisible by $p$, but $\chi_{j}(1)$ is divisible by $|G|_{p},$ this forces $\chi_{j}(y^{-1}) \in J(R),$ while we also have $[G:C_{G}(y)]\chi_{i}(y) \in |G|_{p}R$ (because $$[G:C_{G}(y)]\frac{\chi_{i}(y)}{\chi_{i}(1)}$$ is an algebraic integer). Thus we may ignore class representatives $y$ for which $|C_{G}(y)|$ is divisible by $p$, and deduce that the $(i,j)$-entry of $LM$  is congruent to $\delta_{ij}$ (mod $J(R)$), as claimed.

\medskip
Clearly, then $$(LM)^{p^{n}} \to I_{f_{0} \times f_{0}}$$ as $n \to \infty$ (notice, by the way, that $LM$ has rational entries).
Now the sequence $(ML)^{p^{n}}$ is Cauchy, with idempotent limit $E.$

\medskip
Now $$(ML)^{1+p^{n}} = M(LM)^{p^{n}}L$$ for each $n,$ so the sequence $[(ML)^{1+p^{n}}]$ converges to limit $M I_{f_{0} \times f_{0}}L = ML.$But this sequence also has limits $MLE$ and $EML.$ Thus $E$ is the identity element of the subalgebra of $M_{\ell_{0} \times \ell_{0}}(R)$ generated by $ML.$

\medskip
The number $\ell_{0} = \ell_{0}(G)$ may also be interpreted as the $p$-rank of the Cartan matrix of $RG$, while Alperin's weight conjecture (in its non-blockwise form) was shown in [8] to be equivalent to proving that the number of blocks of defect zero of $RG$ is always equal to 
$$\sum_{ \sigma \in \mathcal{S}_{p}(G)/G)} (-1)^{|\sigma|}\ell_{0}(G_{\sigma}),$$ with an appropriate refinement for the blockwise version.

\medskip
Before leaving this section, to tie this forthcoming material more explicitly to the Green ring formulation of earlier sections, we remind the reader that the number of virtual fixed points on the Lefschetz conjugation module $\mathcal{L}_{c}(RG)$ is equal to 
$$\sum_{ \sigma \in  \mathcal{S}_{p}(G)/G} (-1)^{|\sigma|} \ell_{0}(G_{\sigma}),$$ and that the blockwise version of Alperin's weight conjecture is equivalent to proving that the number of virtual fixed points of $\mathcal{L}_{c}(B)$ is equal to zero for each block $B$ of positive defect of $RG$. 

\medskip
Let us now combine the methodology of the earlier sections with some of the techniques discussed in this section to prove without Clifford theory or the theory of nilpotent blocks that finite groups which have a normal $p$-complement satisfy Alperin's weight conjecture (this result is well-known, more generally, for $p$-solvable groups)

\medskip
\noindent {\bf Lemma 7.4:} \emph{ Suppose that  $G$ has a normal $p$-complement. Then $G$ satisfies the blockwise version of Alperin's weight conjecture.}

\medskip
\noindent {\bf Proof:} Notice that every section of $G$ has a normal $p$-complement. By the results of [8],  it follows that we need to prove that whenever $B$ is a block of positive defect of $RG$ and $\tilde{B}$ is the corresponding block of $\mathbb{FG},$ , then we have $$\sum_{ \sigma \in  \mathcal{S}_{p}(G)/G} (-1)^{|\sigma|} \ell_{0}(B_{\sigma})=0$$ where $$\ell_{0}(B_{\sigma}) = {\rm dim}_{\mathbb{F}}({\rm Tr}_{1}^{G_{\sigma}}(\tilde{B}_{\sigma})) ={\rm dim}_{\mathbb{F}}(1_{\tilde{B}_{\sigma}} {\rm Tr}_{1}^{G_{\sigma}}(\mathbb{F}G_{\sigma})). $$

\medskip
Whenever $\sigma$ is a non-empy chain in $\mathcal{S}_{p}(G),$ we note that (as $G$ has a normal $p$-complement), each $p$-regular element of $G_{\sigma}$ centralizes $O_{p}(G_{\sigma}) \neq 1,$ so that $\ell_{0}(G_{\sigma}) = 0$ and $${\rm Tr}_{1}^{G_{\sigma}}(\mathbb{F}G_{\sigma}) = 0.$$ Hence we have $$\sum_{ \sigma \in  \mathcal{S}_{p}(G)/G} (-1)^{|\sigma|} \ell_{0}(B_{\sigma})= \ell_{0}(B) \geq 0.$$ 

\medskip
The non-blockwise version of this argument shows that $$\sum_{ \sigma \in  \mathcal{S}_{p}(G)/G} (-1)^{|\sigma|} \ell_{0}(G_{\sigma})= \ell_{0}(G).$$

\medskip
If $B^{\prime}$ is any block of defect zero of $RG$, then we know that $$\sum_{ \sigma \in  \mathcal{S}_{p}(G)/G} (-1)^{|\sigma|} \ell_{0}(B^{\prime}_{\sigma})= 1.$$ 

\medskip
We now claim that it suffices to prove that $f_{0}(G) = \ell_{0}(G)$. For if we have $\ell_{0}(G) = f_{0}(G)$ then we must have $\ell_{0}(B) = 0$ for every block $B$ of positive defect of $RG$ and $\ell_{0}(B^{\prime}) = 1$ for every block $B^{\prime}$ of defect zero of $RG$. 

\medskip
But, returning to the matrix  $[a_{ij}]$ defined at the beginning of this section, it is clear that when $y_{i}$ and $y_{j}$ are (necessarily $p$-regular) elements of $p$ defect zero in $G$, we have $\frac{|\Omega_{ij}|}{|S|} = 0$ if $y_{i},y_{j}$ are not conjugate, and we have $$\frac{|\Omega_{ii}|}{|S|} = \frac{[G:C_{G}(y_{i})]}{|S|}$$ if $y_{i} = y_{j}.$ Hence the $\ell_{0} \times \ell_{0}$ matrix $[a_{ij}]$ is diagonal with non-zero entries on the main diagonal, so has rank $\ell_{0}.$ On the other hand, the rank of $[a_{ij}]$ is $f_{0}(G),$ the number of blocks of defect zero of $RG$. Hence we do have the required equality $\ell_{0}(G) = f_{0}(G),$ and the proof of the lemma is complete.

\medskip
In a similar vein, we may prove:

\medskip
\noindent {\bf Lemma 7.5:} \emph{ Suppose that $p = 2$, and that $G$ is a finite group which contains no Klein $4$-subgroup. Then the number of blocks of defect zero of $RG$ is equal to the number of (necessarily $2$-regular) conjugacy classes of $2$-defect zero of $G$, and each block $B$ of positive defect of $RG$ is a block of Lefschetz type.}

\medskip
\noindent {\bf Proof:} It is well-known that a Sylow $2$-subgroup $S$ of $G$ is cyclic, quaternion of order $8$, or generalized quaternion of order at least $16$. By the Brauer-Suzuki theorem in the (generalized) quaternion Sylow $2$-subgroup case, or by transfer in the cyclic Sylow case,  we have $G = O_{2^{\prime}}(G)C_{G}(t)$ for the unique involution $t \in S.$ Note that every block of positive defect of $RG$ has a defect group containing $t$. By the main theorem of [13], there is a bijection between blocks of positive defect of $RG$ and blocks of $RC_{G}(t)$ such that correspnding blocks have the same Cartan matrix. Since the Cartan matrix of every block of $RC_{G}(t)$ has every entry even, we see that $\ell_{0}(B) = 0$ for every block $B$ of positive defect of $RG$. From this it follows that all determinantal divisors $1$ of the Cartan matrix of $RG$ come from blocks of defect zero, and  that the number of blocks of defect zero of $G$ is $\ell_{0}(G).$ 

\medskip
If $B$ is a block of $RG$ of positive defect, then we may choose a defect group $D$ of $B$ which contains $t$. If $\sigma$ is a chain of non-trivial $2$-subgroups of $G$, then we chosse a conjugate chain so that $G_{\sigma} \leq C_{G}(t).$ Then if $B_{\sigma}$ is non-zero, then $\ell_{0}(B_{\sigma}) = 0,$ since every Cartan invariant of $B_{\sigma}$ is even. Hence we have $$\sum_{ \sigma \in \mathcal{S}_{p}(G)/G}(-1)^{|\sigma|}\ell_{0}(B_{\sigma}) = 0,$$ so that $B$ is of Lefschetz type.

\medskip
Returning to the general finite group $G$,  each block of defect zero of $RG$ has Cartan matrix $[1]$, so we need to understand the relationship between $\ell_{0}(G)$ and the number of blocks of defect zero of $RG$, if we hope to verify the alternating sum \emph{conjectural} formula for the number of blocks of defect zero in terms of $\ell_{0}$. This is the subject matter of the next section, where our main concern will be to describe the $p$-rank of the Cartan matrix of each block $B$ of $RG$ via the theory of lower defect groups, as originally developed by R. Brauer, and later refined (separately) by M. Brou\'e, [2] and J. B. Olsson, [9], and also J.A. Green,  [6].

\medskip
It may also be helpful to recall while reading the next section that the virtual character of $G$ afforded by $\mathcal{L}_{c}(RG)$ is 
$$\sum_{\sigma \in \mathcal{S}_{p}(G)/G} (-1)^{|\sigma|} {\rm Ind}_{G_{\sigma}}^{G}[\sum_{ \chi \in {\rm Irr}(G_{\sigma})}(\overline{\chi} \chi )].$$

\medskip
\section{On the multiplicity of $\{1_{G}\}$ as a lower defect group}

\medskip
Here, we examine the multiplicity of $\{1_{G}\}$ as a lower defect group, its relationship to the Cartan matrix, and to the Alperin weight conjecture, (mainly) for blocks of positive defect.

\medskip
The Cartan matrix $C$ of $RG$ has a direct sum decomposition according to the blocks of $RG$. Since all rational primes other than $p$ are invertible in $R$, and since the Brauer characters of the simple $\mathbb{F}G$-modules are linearly independent over $R$, we may conclude from the matrix equation
$$\overline{\Phi}^{t}C \Phi = {\rm diag} \{|C_{G}(y_{i})|: 1 \leq i \leq \ell \}$$ (where $\Phi$ is the Brauer character table of $G$, written as an $\ell \times \ell$ matrix, and $\{ y_{i}: 1 \leq i \leq \ell \}$ is a set of representatives for the conjugacy classes of $p$-regular elements of $G$), that the determinantal divisors of $C$ as an integral matrix are (counting multiplicities) $$\{|C_{G}(y_{i})|_{p} : 1 \leq i \leq \ell \},$$ where we label these classes so that $|C_{G}(y_{i})|_{p}$ divides $|C_{G}(y_{j})|_{p}$ whenever $i < j.$ This comes from the observation that the matrix $\overline{\Phi}^{t}C \Phi $ is already in Smith Normal Form over the principal ideal domain  $R$, and $\mathbb{Z}$ embeds as a subring of $R$ (while the determinantal divisors in both cases must be equal up to multiplication by units in $R$). Notice then that the multiplicity of $1$ as a determinantal divisor of the Cartan matrix $C$ is $\ell_{0}(G),$ the number of $p$-defect zero classes of $G$ (this is also the $p$-rank of $C$).

\medskip 
R. Brauer used this decomposition to define lower defect groups of blocks of $RG$. Each lower defect group is a Sylow $p$-subgroup of the centralizer of some $p$-regular element $y_{i}$ of $G$ (strictly speaking, this refers to the multiplicity of lower defect groups in the $p$-section of the identity). Brauer had a determinant argument to assign $p$-regular conjugacy classes to blocks of $RG$. Later, Brou\'e, [2], and Olsson, [10], (separately) developed different treatments of lower defect groups in terms of ideals of $Z(\mathbb{F}G)$, and we adopt their viewpoints here.

\medskip
For a block $B$ of $\mathbb{F}G$ with defect group $D$, it is a well known consequence of Brauer's Second Main Theorem that whenever $C$ is a class sum of a $p$-regular conjugacy class, the product $1_{B}C$ is an $\mathbb{F}$-linear combination of $p$-regular conjugacy class sums (a similar statement lifts to $Z(RG)$).

\medskip
For each $p$-subgroup $P$ of $G$, the centre of the group algebra, $Z(\mathbb{F}G)$ has an ideal $Tr_{P}^{G}([\mathbb{F}G]^{P}),$ but when dealing with $p$-regular conjugacy classes, we have to content ourselves with dealing with $\mathbb{F}$-vector spaces. The multiplicity of $P$ as a lower defect group (in the $p$-section of the identity) for $\mathbb{F}G$ is 

$${\rm dim}_{\mathbb{F}}\left({\rm Tr}_{P}^{G}([\mathbb{F}G_{p^{\prime}}]^{P}) \right) - {\rm dim}_{\mathbb{F}}\left( \sum_{Q < P}({\rm Tr}_{Q}^{G}([\mathbb{F}G_{p^{\prime}}]^{Q} ) \right).$$ Similarly, for each block $B$ of $\mathbb{F}G$,
the multiplicity of $P$ as a lower defect group (in the $p$-section of the identity) for $B$ is 

$${\rm dim}_{\mathbb{F}}\left({\rm Tr}_{P}^{G}(1_{B}[\mathbb{F}G_{p^{\prime}}]^{P}) \right) - {\rm dim}_{\mathbb{F}}\left( \sum_{Q < P}({\rm Tr}_{Q}^{G}(1_{B}[\mathbb{F}G_{p^{\prime}}]^{Q} ) \right).$$ In particular, the multiplicity of $1_{G}$ as a lower defect group for $B$, denoted by $\ell_{0}(B),$ simplifies to  $${\rm dim}_{\mathbb{F}}\left({\rm Tr}_{1}^{G}(\tilde{B})\right), $$ where $\tilde{B}$ is the block of $\mathbb{F}G$ corresponding to $B$. This definition may be extended to sums of blocks.

\medskip
By the results of [8], Alperin's weight conjecture is equivalent to proving that $$\sum_{\sigma \in \mathcal{ S}_{p}(G)/G} (-1)^{|\sigma|}\ell_{0}(B_{\sigma})  = 0$$ for every block $B$ of $RG$ of positive defect (for every finite group $G$). This now becomes a  question about the multiplicity of $1$ as a lower defect group for blocks $B$ of positive defect and their Brauer correspondents. We note also that $\ell_{0}(B)$ is the multiplicity of the projective cover of the trivial module as a summand of $B$ viewed as $RG$-module under the conjugation action of $G$.

\medskip
Now let us work with $RG$ for a while. The next Lemma is due in essence to Osima.

\medskip
\noindent {\bf Lemma 8.1:} \emph{ Let $B$ be a sum of blocks of $RG$, and let $K_{y}$ denote the class sum in $RG$ of the conjugacy class containing the $p$-regular element $y \in G$. Then $K_{y}1_{B}$ is an $R$-linear combination of class sums of $p$-regular elements, and for each $p$-regular element $z \in G,$  the coefficient of $z$ in $K_{y}1_{B}$
is $$\frac{1}{|C_{G}(y)|}\left( \sum_{\chi \in {\rm Irr}(B)}\chi(y)\chi(z^{-1}) \right).$$ In particular, this last expression lies in $R$. Furthermore, if a Sylow $p$-subgroup of $C_{G}(z)$ is not contained in a conjugate of $C_{G}(y)$, then this last expression lies in $pR \subseteq J(R).$}

\medskip
\noindent {\bf Proof:} Note that in $Z(\mathbb{K}G),$ we have $$K_{y}1_{B} = \sum_{\chi \in {\rm Irr}(B)}\frac{[G:C_{G}(y)]\chi(y)}{\chi(1)} \left( \frac{\chi(1)}{|G|}\right) \left( \sum_{g \in G}\chi(g^{-1})g \right),$$ so, whether or not $z \in G$ is $p$-regular, the coefficient of $z$ in $K_{y}1_{B}$ is as stated. However, by Brauer's Second Main Theorem and orthogonality relations derived from it, we know that $$\sum_{\chi \in {\rm Irr}(B)} \chi(y)\chi(z^{-1})= 0$$ whenever $z \in G$ is $p$-singular.

\medskip
For the last claim, note that if the coefficient of $z$ in $K_{y}1_{B}$ lies outside $pR$, then there must be a class sum $K_{x}$ appearing with non-zero coefficient in $1_{B}$ such that the coefficient of $z$  in $K_{y}K_{x}$ is not divisible by $p$. However, $C_{G}(z)$ acts by conjugation on ordered pairs $(y^{\prime},x^{\prime})  \in y^{G} \times x^{G}$ with $y^{\prime}x^{\prime} = z,$ so some orbit must have length prime to $p$, and, in particular, a Sylow $p$-subgroup of $C_{G}(z)$ must centralize some conjugate of $y$, so the result follows.

\medskip
Now let $B$ be a direct sum of blocks of $RG$, and let $\tilde{B}$ denote the corresponding sum of blocks of 
$\mathbb{F}G.$ Let $\{y_{i}: 1 \leq i \leq \ell \}$ be a set of representatives for the $p$-regular conjugacy classes of $G$, and let $K_{i}$ denote the class sum of $y_{i}$ in $Z(RG)$ for each $i$, and $\tilde{K_{i}}$ denote the class sum of $y_{i}$ in $\mathbb{F}G.$ Let $S_{{\rm reg}}$ denote the $R$-span of $\{K_{i}\}$ in $Z(RG)$, and let $\tilde{S}_{\rm reg}$ denote its image in $Z(\mathbb{F}G).$

\medskip
Note that $X \mapsto 1_{B}X$ is an idempotent $R$-module homomorphism from $S_{{\rm reg}}$ to $Z(B) \cap S_{{\rm reg}}$  and if we let $A(B)$ denote its $\ell \times \ell$-matrix with respect to the given $R$-basis, the $(i,j)$-entry of $A(B)$ is $$\sum_{\chi \in {\rm Irr}(B)} \frac{\chi(y_{i}) \chi(y_{j}^{-1})}{|C_{G}(y_{i})|}$$ by Lemma 8.1. Since $A(B)$ is an idempotent matrix, its rank is its trace, and this is $\ell(B)$ (or, more precisely $\ell(\tilde{B})$), as already noted in [8], (in the single block case). Its image, $A(\tilde{B}),$  a matrix with entries in $\mathbb{F},$ is also idempotent of rank $\ell(B)$.

\medskip
When $B$ is a single block of $RG$ and $B^{\prime}$ is a different block of $RG$, then we have $$A(B)A(B^{\prime}) = A(B^{\prime})A(B) = [0].$$
It follows that $I_{\ell \times \ell}  = \sum_{B}A(B)$ as a sum of mutually orthogonal idempotent matrices, were $B$ runs over the blocks of $RG$.

\medskip
Working over $\mathbb{F}$, we can perform a similar argument with the subspaces $${\rm Tr}_{Q}^{G}((\mathbb{F}G)_{p^{\prime}}^{Q})$$ of $Z(\mathbb{F}G)$, for $p$-subgroups $Q$ of $G$. We illustrate with the case $Q= \{1_{G}\}$ for the sake of exposition.

\medskip
Let us label the $p$-regular classes so that $|C_{G}(y_{i})|$ is prime to $p$ for $1 \leq i \leq \ell_{0},$ and let $S_{0}$ denote the $\mathbb{F}$-span
of $\{ \tilde{K_{i}} : 1 \leq i \leq \ell_{0} \}.$ Then for each block $\tilde{B}$ of $\mathbb{F}G$, the linear map $X \mapsto 1_{\tilde{B}}X$ is an idempotent linear mapping, and its matrix $A_{0}(\tilde{B})$ with respect to the given  basis is an idempotent matrix whose rank is the multiplicity of $1_{G}$ as a lower defect group for $\tilde{B},$ which is also the $p$-rank of the Cartan matrix for $B$ (that is to say, the number of determinantal divisors $1$ of the Cartan matrix of $B$).

\medskip
For a suitable power of $p$, say $q$, the limit $$\lim_{n \to \infty} (A_{0}^{\prime}(B))^{q^{n}}$$ exists, and is an idempotent matrix whose trace is the rank of $A_{0}(B),$ for any preimage $A_{0}^{\prime}(B)$ of $A_{0}(B)$ in $M_{\ell_{0} \times \ell_{0}}(R).$  By this lifting process, we may write $$I_{\ell_{0} \times \ell_{0}} = \sum_{B} \tilde{A_{0}}(B),$$ where the  $\tilde{A_{0}}(B)$ are mutually orthogonal idempotent matrices whose rank (and trace) is the multiplicity of one as a lower defect group for the block $B$ of $RG$. In practice, it is usually sufficient to deal with the matrices $A_{0}(B)$ in 
$M_{\ell_{0} \times \ell_{0}}(\mathbb{F}),$ as these have the correct rank, though their traces lie in $\mathbb{F}$.

\medskip
We note that one key difference between blocks of defect zero and blocks of positive defect is that when $B$ is a block of defect zero, ${\rm Tr}_{1}^{G}(B)$ is not contained in  $J(Z(B)),$ while if $B$ has positive defect, all elements of ${\rm Tr}_{1}^{G}(B)$ lie in $J(Z(B)).$ Unlike the situation of Section 7, for a block $B$ of positive defect of $RG$, the fact that $\ell_{0}(B)$ is the rank of an idempotent matrix with entries in $\mathbb{F}$ has nothing to do with idempotent elements of ${\rm Tr}_{1}^{G}(B)$, since there are no such elements..

\medskip
\noindent {\bf Remark 8.2:} Let $G = A_{5}$ and $p = 3.$ Let $B$ be the principal block of $RG$. There are three $3$-regular classes of $G$ whose representatives have centralizers of order prime to $3$: these are the classes of $y_{1} = (12)(34)$, as well as $y_{2} = (12345)$ and $y_{3} = (13524)$.

\medskip
From the character table of $G$, we see that the matrix $A_{0}(B)$ is equal  (over $\mathbb{F}$) to $$\left(\begin{array}{clcrc} 1/2 & 1/4& 1/4\\ 1/5 & 2/5 & 2/5\\ 1/5 & 2/5 & 2/5 \end{array} \right),$$ which makes sense. Note that this matrix is not idempotent if read over $R$, since its trace is not even an integer. However, a preimage of $A_{0}(B)$ with entries in $R$ is $$A_{0}^{\prime}(B) = \left(\begin{array}{clcrc} -1 & 1& 1\\ -1 & 1 & 1\\ -1 & 1 & 1 \end{array} \right),$$ and, somewhat serendipitously, this matrix is already idempotent, of rank one and trace one, so we may conclude that the multiplicity of $1$ as a lower defect group for the block $B$ is one, which is correct.

\medskip
\noindent {\bf Acknowledgement:} We are grateful to Burkhard  K\"ulshammer for pointing out an inaccuracy in an earlier version of this work, and to Gunter Malle for a number of typographical corrections

\newpage
\medskip
\begin{center} {\bf References} \end{center}

\medskip
\noindent [1]  Alperin, J.L., \emph{Weights for finite groups},  {\bf in}: The Arcata Conference on Representations of Finite Groups (\emph{ed. P. Fong}), AMS Proceedings of Symposia in Pure Mathematics, {\bf 47}, (1987).

\medskip
\noindent [2] Brou\'e, M., \emph{Les $l$-blocs des groupes finis et leurs invariants}, Inventiones Mathematicae, {\bf 56}, (1980), 117-136.

\medskip
\noindent [3] Conway, J.H., Curtis, R.T., Norton, S.P., Parker, R.A., \emph{Atlas of Finite Groups}, Clarendon Press, Oxford, (1985).

\medskip
\noindent [4] Fong, P., \emph{On the characters of $p$-solvable groups}, Transactions of the American Mathematical Society, {\bf 98}, (1961),263-284.

\medskip
\noindent [5] Green, J.A., \emph{Some remarks on defect groups}, Mathematische Zeitschrift, {\bf 107},(1968),133-150.

\medskip
\noindent [6] Green, J.A., \emph{Multiplicities, Scott Modules and lower defect groups}, Journal of the London Mathematical Society,{\bf 28},(1983),282-292.

\medskip
\noindent [7] Harris,M.E., \emph{Blocks of small defect}, Journal of Algebra, {\bf 104}, (1986), 167-195.

\medskip
\noindent [8] Kn\"orr, R. and Robinson,G.R., \emph{Some remarks on a conjecture of Alperin}, Journal of the London Mathematical  Society,{\bf 39},(1989),48-60.

\medskip
\noindent [9] K\"ulshammer, B., \emph{Bemerkungen \"uber die Gruppenalgebra als symmetrische Algebra, III}, Journal of Algebra, {\bf 88},(1984),279-291.

\medskip
\noindent [10] Olsson, J.B., \emph{Lower defect groups in blocks}, Communications in Algebra,{\bf 8},(1980),261-288.

\medskip
\noindent [11] Reynolds, W.F., \emph{Blocks and normal subgroups of finite groups}, Nagoya Mathematical Journal, {\bf 22}, (1963), 15-32.

\medskip
\noindent [12] Robinson,G.R., \emph{The number of blocks with a given defect group}, Journal of Algebra,{\bf 84}, (1983), 493-502.

\medskip
\noindent [13] Robinson,G.R.,\emph{On Cartan matrices of finite groups with a certain factorization},  Representation Theory II,Groups and Orders,  Springer Lecture Notes 1178,152-176, Berlin (1986).

\medskip
\noindent [14] Robinson, G.R., \emph{Some remarks on permutation modules}, Journal of Algebra,{\bf 118},(1988), 46-62.

\medskip 
\noindent [15] Robinson, G.R., \emph{A note on Alperin's conjecture}, Journal of Algebra,{\bf 165},(1994),602-605.

\medskip 
\noindent [16] Robinson, G.R., \emph{A generalized character related to the $p$-local structure and representation theory of a finite group},
Journal of Algebra, {\bf 690},(2026),37-74.

\medskip
\noindent [17] Webb, P.J., \emph{Subgroup complexes}, {\bf in}: The Arcata Conference on Representations of Finite Groups (\emph{ed. P. Fong}), AMS Proceedings of Symposia in Pure Mathematics, {\bf 47}, (1987).

\end{document}